\documentclass[smallcondensed]{svjour3}  
\usepackage{amsmath,color,amssymb,graphicx,hyperref,lscape,bbm}
\definecolor{darkblue}{rgb}{0.0,0.0,0.3}
\hypersetup{colorlinks,breaklinks,
            linkcolor=darkblue,urlcolor=darkblue,
            anchorcolor=darkblue,citecolor=darkblue}

\begin{document}
\newcommand{\R}{\mathbb{R}}
\newcommand{\N}{\mathbb{N}}
\newcommand{\Z}{\mathbb{Z}}
\newcommand{\I}{\mathcal{I}}
\newcommand{\LL}{\mathcal{L}}
\newcommand{\F}{\mathcal{F}}
\newcommand{\x}{\mathbf{x}}
\newcommand{\X}{\mathbf{X}}
\newcommand{\xxi}{\mathbf{\xi}}
\newcommand{\y}{\mathbf{y}}
\newcommand{\m}{\mathfrak{m}}
\newcommand{\vv}{\mathbf{v}}
\newcommand{\e}{\mathbf{e}}
\newcommand{\ud}{\mathrm{d}}
\newcommand{\hra}{\hookrightarrow}
\bibliographystyle{plain}

\title{A nonsmooth program for jamming hard spheres}
\author{Peter Hinow}
\institute{Department of
Mathematical Sciences, University of Wisconsin -- Milwaukee, P.O.~Box 413,
Milwaukee, WI 53201, USA; \email{hinow@uwm.edu}}

\date{\today}
\maketitle
\begin{abstract}
\begin{sloppypar}
We study packings of $n$ hard spheres of equal radius in the $d$-dimensional unit cube. We present a nonsmooth function whose local extrema are the radii of jammed packings (where no subset of spheres can be moved keeping all others fixed) and show that for a fixed number of spheres there are only finitely many radii of such jammed configurations. We propose an algorithm for the maximization of this maximal radius function and present examples for 5 - 8 disks in the unit square and 4 - 6 spheres in the unit cube. The method allows straightforward generalization to packings of spheres in other compact containers.
\end{sloppypar}

\end{abstract}
\section{Introduction}\label{Intro}
Originally introduced as models for monatomic liquids, granular media and
glasses, random dense sphere packings have been of high interest to experimental
and theoretical physicists as well as mathematicians for more than 50 years. The classical experiments of
Bernal \cite{Bernal1960a,Bernal1960b} and Scott \cite{Scott1960} from 1960 were
to pour small balls of equal radius into a container and shake the container for
the packing to densify. Although this is a parameter-dependent protocol (the
speed of pouring and the shaking amplitude and frequency), it was
observed that there exists a ``universal'', highly reproducible packing fraction
of $\approx 0.63-0.64$ for disordered monodisperse spheres \cite{Kamien}. Later,
Lubachevsky and Stillinger \cite{Lubachevsky1990} introduced the ``stochastic
billiard'', a computational tool to create dense sphere packings. There, a
fixed number of spheres move within a container while their radii increase at
a rate common to all. The spheres collide elastically with each other and with the
walls of the container if such are present. At every collision, the kinetic
energy of the collision partners is increased to ensure their separation. This
procedure is stopped once either the collision frequency or the pressure (the
sum of the squared velocities) exceed a certain threshold. This is also a
parameter dependent protocol, in that the radius growth rate has to be set by
the user. Very slow growth rates result in higher packing fractions. However,
the spheres are then arranged in local patches of crystalline, and thus regular substructures.
Applications of random dense sphere packings are manifold, from modeling the
burning of solid rocket propellant \cite{Knott}, predicting the drug release
kinetics of matrix tablets in pharmaceutical science \cite{DCDS_B} to phase transitions in statistical mechanics
\cite{Diaconis09,Loewen}, to name only a few.

Torquato \textit{et al.}~\cite{Torquato} argued that the concept of random close
packing of spheres is ill defined. Clearly, the properties of density and
randomness are at variance with each other and an increase in one must come at
the cost of a decrease of the other. The authors of \cite{Torquato} propose to
balance an increase in the packing fraction $\phi$ with an increase of a
suitably chosen ``order parameter'' $\psi$. They suggest that there are
``maximally random jammed'' (MRJ) states, namely jammed sphere packings
that minimize that order parameter. Following Torquato \textit{et
al.}~\cite{Torquato,Donev2004a}, a packing of finitely many spheres
in a container is \textit{(collectively) jammed} if no subfamily of spheres can
be displaced continuously while fixing the positions of all other spheres. The present paper is motivated by the question for which radii there are jammed configurations of $n$ spheres in the $d$-dimensional unit cube. Many authors have studied densest packings of  50 - 100 spheres in  squares, triangles, disks and the cube, for certain numbers of spheres also with computer-aided proofs of optimality, see  for example   \cite{Boll,Croft,Gensane,Joos,Maranas,Nurmela97,Nurmela99,Peikert}  and the references therein. Nurmela and \"{O}stergaard \cite{Nurmela97} achieved this by minimizing the (smooth) discrete Riesz \mbox{$s$-energy} of a configuration $\x$ of sphere centers $\x_i$ (for more details on the notation see definition \eqref{defx} below)
\begin{equation*}
E_s(\x) = \sum_{i<j} |\x_i-\x_j|^{-s}
\end{equation*}
in the limit $s\to\infty$, when this energy approaches the nearest neighbor distances \cite{Hardin}, see also the recent work by Carlsson \textit{et al.}~\cite{Carlsson}. Rather than seeking only densest packings, here we emphasize finding as many jammed packings as possible and work directly with a nonsmooth  function to be optimized.

In this paper we take up the idea from \cite{Peikert} of the maximal radius function $G$ (see Equation \eqref{lower_envelope} below)  for finite configurations and show that jammed packings correspond to its local maxima. We prove that for every finite number of spheres in $[0,1]^d$ there are only finitely many radii at which jamming can occur.  Then we apply techniques from nonsmooth analysis \cite{Clarke1983,Makela92} to optimize $G$. Our algorithm is based on a line search and has similarities to the gradient sampling algorithm by Burke {\it et al.}~\cite{Burke}. To find the direction of the line search we exploit the special structure of the objective function whose gradients are readily calculated. This optimization procedure is repeated many times with random starting points. With this parameter-free approach we are able to find many jammed packings for up to eight disks in the square and six spheres in the cube together with the minimal polynomials for the corresponding radii (thereby proving their existence). To the best of our knowledge, many of these configurations have not been reported elsewhere in the literature. Our method can easily be generalized to sphere packings in arbitrary convex polytopes, within a sphere or to spherical caps on a sphere. The paper ends with a discussion of open problems for future research.

\section{Jammed sphere packings and the maximal radius function}\label{Theorem}

Throughout the paper, $|\,\cdot\,|$ denotes the Eucliden norm. Let the number of spheres $n\ge 2$ be fixed and let
\begin{equation}\label{defx}
\x=(\x_1,\,\dots\,,\x_n)=
(x_{11},\,\dots\,,x_{1d},\,\dots\,,x_{n1},\,\dots\,,x_{nd})
\end{equation}
be a collection of $n$ points $\x_j\in\R^d$, $d\ge2$. We say that $\x$ is
$r$-\textit{admissible} for $r\ge 0$ if 
$|\x_i-\x_j| \ge 2r$
for all $i,j=1,\dots,n$ with $i\neq j$ and 
$r\le x_{ik} \le 1-r$
for all $i=1,\dots,n$ and $k=1,\dots,d$. The latter conditions can be easily generalized to convex polytopes bounded by half-spaces with known normal vectors. Notice that for any such point $\x$, a permutation of the entries $\x_i$ gives another point $\x'$ that determines the same geometric configuration of spheres. In an even stronger sense, one can identify configurations that are obtained from each other by an application of a symmetry of the cube. 
 
Denote by $M_n^d(r)\subset\R^{dn}$ the set of all $r$-admissible
configurations of $n$ distinct points in $[0,1]^d$. When $r=0$, then
$M_n^d(0)$ is path-connected, since any two components $\x_i$ and $\x_j$
can exchange their places along continuous curves, each avoiding all other
particles and the faces of the cube. If $r_1<r_2$, then $M_n^d(r_2)\subset
M_n^d(r_1)$. Moreover, for every $d$ there is a sequence $r_n$ with
$\displaystyle\lim_{n\to\infty} r_n = 0$ such that
$M_n^d\left(r_n\right) = \varnothing$ (large numbers of spheres
imply small radii). As $r$ increases, path-connected components split off (i.e.~any two points in the same component can be joined by a continuous curve)
at certain times and disappear at others and so the topology of $M_n^d(r)$ changes.
We say that $r$ is \textit{critical}, if there exists a $\varepsilon>0$ such
that for the numbers of connected components (denoted by $\beta_0$) we have $\beta_0(M_n^d(s))\neq
\beta_0(M_n^d(t))$ for every $r-\varepsilon<s<r<t<r+\varepsilon$. In
particular, if $\beta_0(M_n^d(s))> \beta_0(M_n^d(t))$, we speak of a
\textit{disappearance} and we call the set
\begin{equation*}
C(r)=M_n^d(r)\setminus \overline{\bigcap_{t>r}M_n^d(t)} 
\end{equation*}
a \textit{critical set} (this is the union of all disappearing connected
components at the critical radius $r$).  

\begin{figure}[ht]
\begin{center}
\includegraphics[width=120mm,height=82mm]{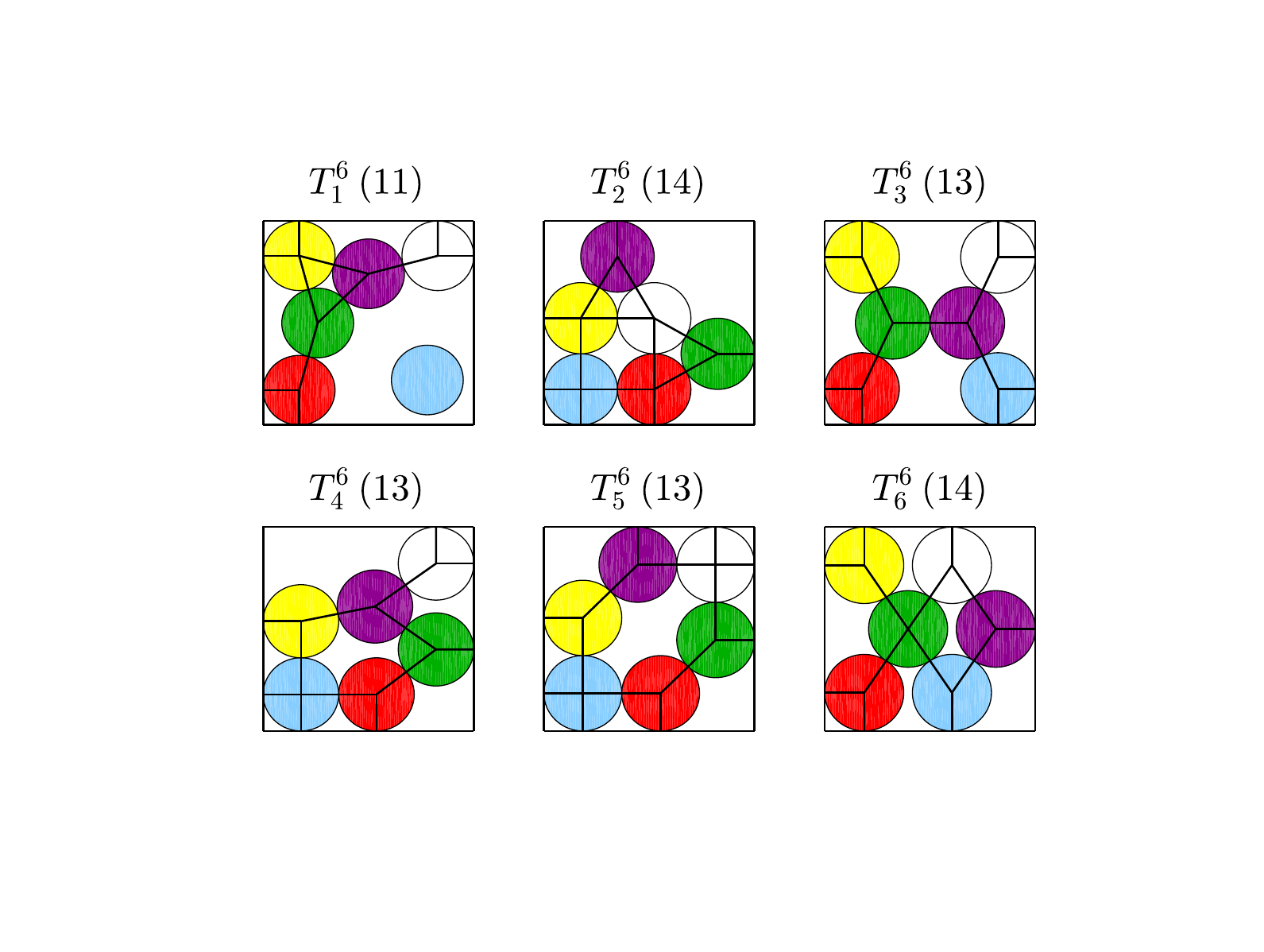}
\caption{
Partially and fully jammed configurations of six
disks in a square with increasing radii from left to right and top to bottom.  
Configuration $T_2^6$ allows another disk of the same size to be added.
}\label{n6}
\end{center}
\end{figure}

\begin{figure}[ht]
\begin{center}
\includegraphics[width=120mm,height=100mm]{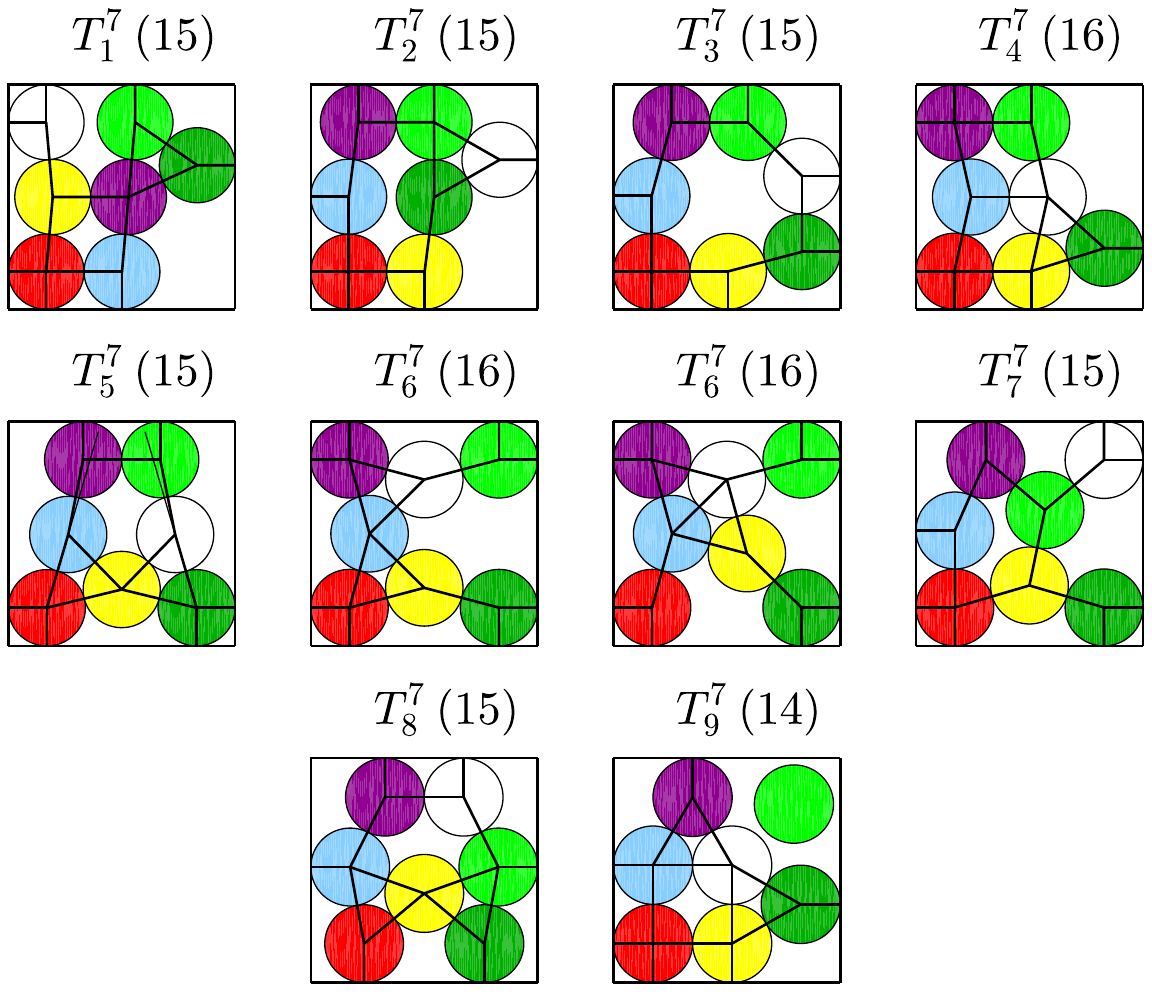}
\caption{Partially and fully jammed configurations of seven disks and splitting events with increasing radii from left to right and top to bottom. The thin solid lines in configuration $T^7_5$ indicate that it is indeed jammed. Configuration $T_6^7$ also exists with the two disks as rattlers. 
}\label{n7}
\end{center}
\end{figure}

For $r > 0$, $M^d_n(r)$ is clearly closed. If $\x \in int\,M^d_n(r)$ and $r>
0$, then no two spheres centered at the points in $\x$ touch and there exists a
$\varepsilon> 0$ such that $\x \in int\,M^d_n(r')$ for every $r'\in [r, r +
\varepsilon)$. Thus $C(r)$ at a critical radius has empty interior. We call a point  $\x\in
C(r)$ a \textit{partially jammed} configuration of spheres and an isolated point
of $C(r)$ a \textit{fully jammed} configuration of spheres. The fact that
$int\,C(r)=\varnothing$ excludes that the spheres can be displaced continuously
so that eventually they all loose contact with the walls and each other and the
radius can be increased. However, it is possible that a
partially jammed configuration leaves room for an unconstrained sphere (a
``rattler'', see Figures \ref{n6} and \ref{n7})\footnote{Figures and Tables are collected at the end of this manuscript.}.  Clearly,  there is always a largest radius $r$ for which $M_n^d(r)$ is not empty.
 
Jammed configurations can be characterized as local maxima of
the maximal radius function that we now introduce. With the notation as in \eqref{defx}, let 
\begin{equation*}
\begin{aligned}
\varphi_{ij}(\x)&= \frac{|\x_i-\x_j|}{2},\quad i=1,\dots,n-1,\,
j=i+1,\dots,n, \\
\psi_{ik}(\x)&= x_{ik},\quad \psi^{ik}(\x)= 1- x_{ik}, \quad i=1,\dots,n,\,
k=1,\dots,d, 
\end{aligned}
\end{equation*}
\begin{sloppypar}\noindent
and denote the set of these $N(n,d):=\frac{n(n-1)}{2}+2nd$ functions by
$\F$.
Define \mbox{$G:[0,1]^{nd}\to[0,\infty)$} by the lower envelope
\end{sloppypar}
\begin{equation}\label{lower_envelope}
G(\x) = \min_{f\in\F}\left\{f(\x)\right\}.
\end{equation}
\begin{sloppypar}\noindent
This is the maximal $r$ such that a sphere of radius $r$ can be
centered at every entry of the $n$-tuple $\x=(\x_1,\,\dots\,,\x_n)$ without
violating any of the other spheres or the walls of the unit cube.
Throughout, we fix an index set \mbox{$\LL=\{1,2,\dots,N(n,d)\}$} and relabel the functions so that $\F=\{f_l\::\:l\in\LL\}$. 
The gradients of the component functions are
\end{sloppypar}
\begin{equation*}
\nabla\phi_{ij}(\x)= \frac{1}{2|\x_i-\x_j|}\left(\dots, 0, \dots,
\underbrace{x_{ik}-x_{jk}}_{k=1,\dots,d},\dots, 0, \dots,
\underbrace{x_{jk}-x_{ik}}_{k=1,\dots,d},\dots, 0,\dots
\right), 
\end{equation*}
where the first stretch begins at the entry $(i-1)d$ and the second
begins at the entry $(j-1)d$.  The gradients of the $\phi_{ij}$ exist wherever $\x_i\neq \x_j$ and are non-zero there.  Trivially,
\begin{equation*}
\nabla\psi_{ik}(\x)= \delta_{((i-1)d+k,l)}, \quad
\nabla\psi^{ik}(\x)= -\delta_{((i-1)d+k,l)},
\end{equation*}
where $\delta_{k,l}$ denotes the Kronecker symbol (1 if $k=l$ and 0 otherwise).

\begin{lemma} If $\x^*$ is a local maximum of $G$, then $\x^*$ is a
partially jammed configuration.  If $\x^*$ is a strict local maximum of $G$,
then it is fully jammed.
\end{lemma}
\textbf{Proof.} We have that $M_n^d(r)=G^{-1}([r,\infty))$. If there
exists a $\delta>0$ such that $G(\x)\le G(\x^*)$ for all $\x\neq\x^*$
with $|\x-\x^*|<\delta$, then none of these $\x$
lies in $M_n^d(G(\x^*)+\varepsilon)$ for any $\varepsilon>0$.  \hfill
$\Box$

We use methods from nonsmooth optimization to find the maxima of $G$, see
\cite{Clarke1983,Makela92} for references. For a
$K$-Lipschitz continuous function, the \textit{generalized directional
derivative} at  $x$ in the direction $v$ is 
\begin{equation*}
f^\circ(\x;\vv) = \limsup_{\y\to \x,t\searrow0}\frac{f(\y+t\vv)-f(\y)}{t};
\end{equation*}
this is bounded from above by $K|\vv|$. The \textit{generalized gradient} or the
\textit{Clarke subdifferential} is the nonempty set
\begin{equation*}
\partial f(\x) = \left\{\zeta\in\R^m\::\:f^\circ(\x;\vv)\ge\zeta\cdot \vv
\:\textrm{for all} \:\vv\in\R^m\right\}.
\end{equation*}
If $f$ has an extremum at $\x$, then $0\in\partial f(\x)$, \cite[Theorem 3.2.5]{Makela92}. For a function $G$ defined by a minimum selection of functions indexed by a set $\LL$ as in Equation \eqref{lower_envelope} let 
\begin{equation}\label{active}
\LL(\x) = \{l\in\LL\::\:f_l(\x)=G(\x)\} 
\end{equation}
be the indices of the functions that realize the common minimum; these are called the \textit{active indices}.
Then the subdifferential of $G$ at a point $\x$ is
\begin{equation*}
\partial G(\x) = \text{conv}\{\nabla f_l(\x) \::\:l\in\LL(\x)\},
\end{equation*}
the (closed) convex hull of the active gradients, this follows for example from \cite[Theorem 2.5.1]{Clarke1983}. Moreover, the generalized directional derivative is 
\begin{equation*}
G^\circ(\x;\vv) = \min\{\nabla f_l(\x)\cdot \vv \::\:l\in\LL(\x)\}.
\end{equation*}
\begin{proposition} For every number of spheres $n$ and dimension $d$, the function $G$ has finitely many local maxima.
\end{proposition}
\textbf{Proof.} This follows from the fact that every subset of $\F$ can be active at one local maximum at most.  Assume to the contrary that there are local maxima $\x$ and $\x'$ with $\LL(\x)<\LL(\x')$ and $G(\x)<G(\x')$. As all functions participating in $\F$ are convex, the sublevel sets $f_l^{-1}([0,G(\x)])$ are convex. Since the higher point $\x'$ is a local maximum and hence $0\in \text{conv}\{\nabla f_l(\x') \::\:l\in\LL(\x')\}$, there must be at least one gradient that  satisfies $\nabla f_{l^*}(\x')\cdot(\x'-\x)\le 0$. Because of the convexity of the sublevel sets of $f_{l^*}$ this would imply $G(\x)\ge G(\x')$, a contradiction. 
\hfill $\Box$

In the following we describe a iterative procedure to find local maxima of the maximal radius function $G$. An initial choice $\x^0$ is selected for example randomly from a  uniform distribution in $[0,1]^{nd}$. Given the iteration $\x^k$, let $J(\x^k)$ be the matrix whose rows are the active gradients $\nabla f_l(\x^k)$ for all $l\in\LL(\x^k)$. We seek a
direction $\xxi^k$ that solves the minimization problem 
\begin{equation}\label{find_incr} 
\xxi^k\in\underset{\xxi}{\operatorname{argmin}}\,|J(\x^k) \xxi -\mathbbm{1}|^2,
\end{equation}
where $\mathbbm{1}$ is a vector with $|\LL(\x^k)|$ entries $1$. This implies that the active functions all increase infinitesimally at the same
rate. With such a $\xxi^k$ in hand, we begin a line search in that direction. A
triple $0\le t_1<t_2<t_3$ is a \textit{bracket} of a directional maximum, if
\begin{equation}\label{bracket}
G(\x^k+t_1\xxi^k)\le G(\x^k+t_2\xxi^k)\ge G(\x^k+t_3\xxi^k),   
\end{equation}
and at least one of the inequalities is strict. If the width of the bracket $t_3-t_1$ is sufficiently small, then we set $\x^{k+1}=\x^k+t_2\xxi^k$.  

\begin{figure}[ht]
\begin{center}
\includegraphics[width=70mm,height=35mm]{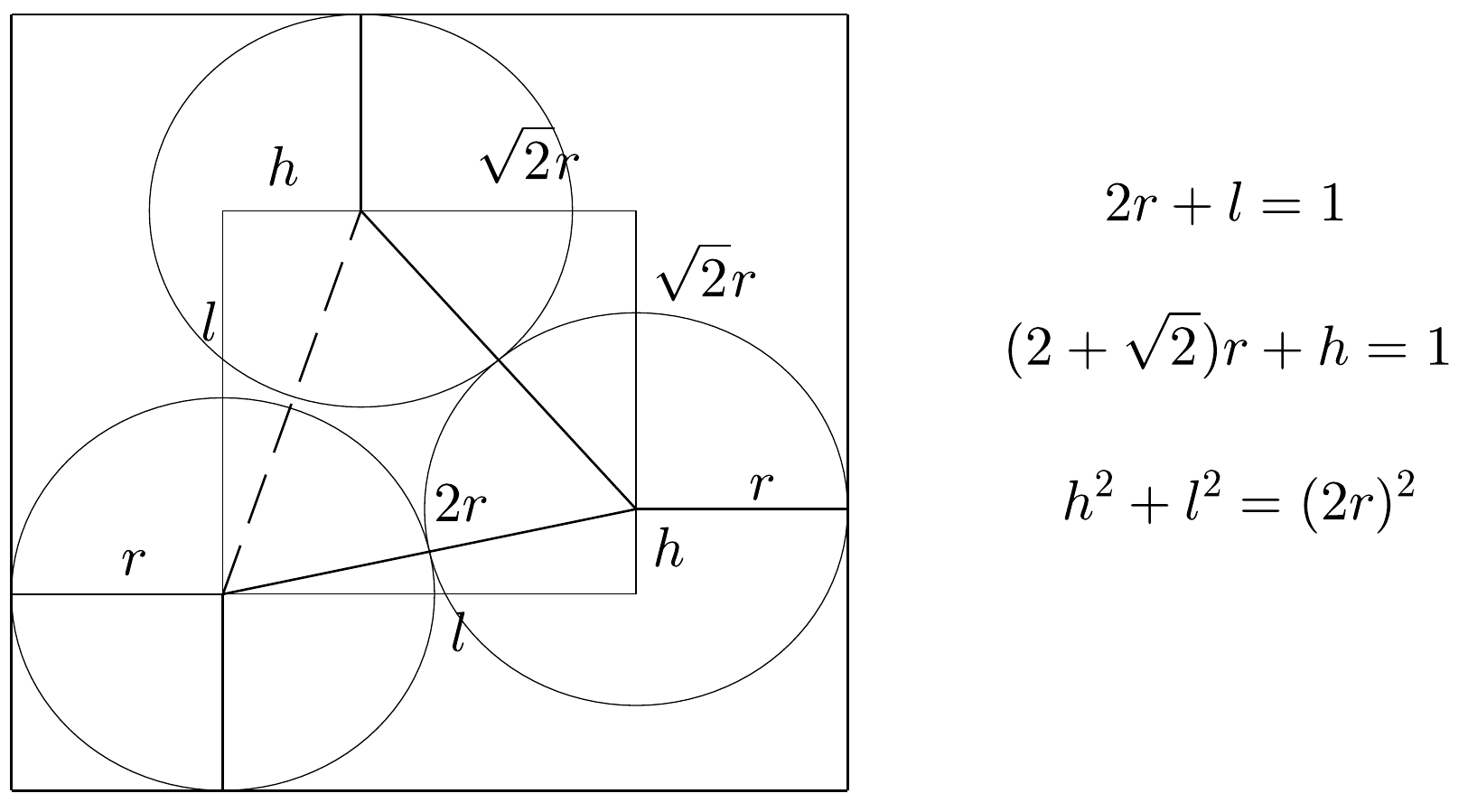}
\caption{From the approximately locally maximal configuration, it is possible to construct the expected contact graph by adding the dashed line. This leads to the system of equations on the right. Solving this system with \textsc{Mathematica} gives $r=\frac{4+\sqrt{2}-\sqrt{6}}{2 \left(3+2 \sqrt{2}\right)}\approx0.254333095$ as the maximal value.}\label{refine}
\end{center}
\end{figure}

There are three issues that need to be addressed in the numerical implementation of this algorithm. First, due to rounding errors we need to relax the equality requirement in the definition \eqref{active} of the set of active functions $\LL$. This is done by replacing \eqref{active} by
\begin{equation*}
\LL_\varepsilon(\x) = \{l\in\LL\::\:f_l(\x)\le G(\x)+\varepsilon\}, 
\end{equation*}
where the tolerance $\varepsilon$ can be adjusted if need be. For example, if active functions are lost after an iteration, then it is helpful to increase $\varepsilon$.
The second is a detection and treatment of saddle points. At the beginning of the line search procedure we set $t_1=0$ and $t_2$ to the machine precision. If instead of the first inequality in \eqref{bracket} we have the reverse strict inequality, then we perturb $\x^k$ by a  random vector whose components are normally distributed with standard deviation $10^{-3}$ and recalculate the direction of increase $\xxi^k$ accordingly.
Finally, the iteration terminates if the optimal function value from the minimization problem \eqref{find_incr} exceeds a certain threshold to be set by the user or a maximum number of iterations has been reached. Once an approximate local maximum has been found, it can be refined by solving a system of linear-quadratic equations determined by the expected contact graph (using e.g.~\textsc{Mathematica}). See Figure \ref{refine} for an example.

\section{Examples}\label{Results}
The algorithm described in Section \ref{Theorem} has been implemented in the open source package \textsc{scilab} \cite{scilab} and together with a \textsc{Mathematica} notebook will be made available on the author's website.

\begin{figure}[ht]
\begin{center}
\includegraphics[width=120mm,height=82mm]{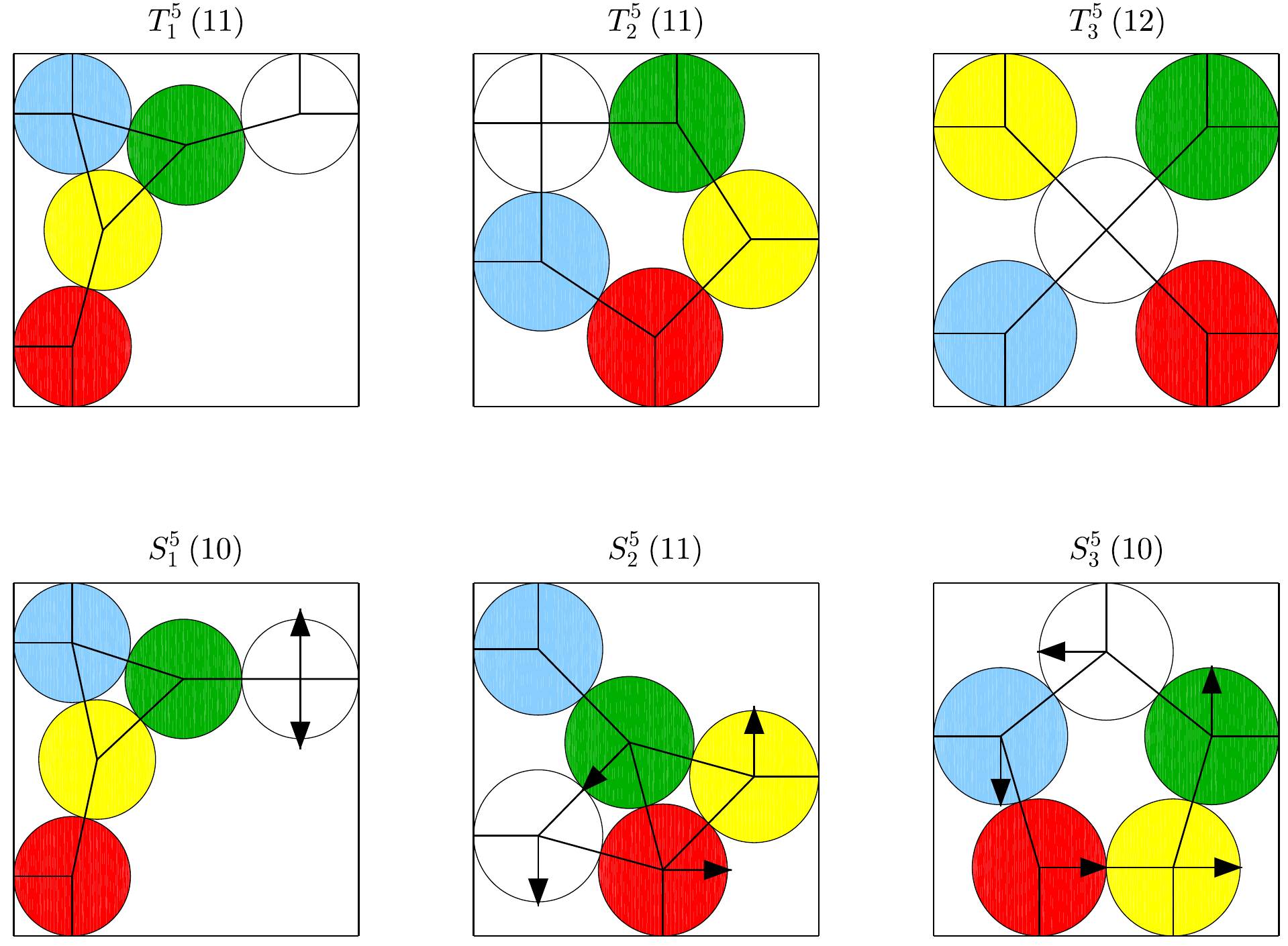}
\caption{\textit{(Top row)} Five fully jammed disks in a square (with increasing
radii from left to right). The straight lines
indicate the contact graphs. The  configuration $T_1^5$ allows up to three more
disks of the same size to be added. Numbers in parentheses indicate the number of functions that realize the common radius.  \textit{(Bottom row)} Critical points
at which connected components of $M_5^2(r)$ split off. The numerical ordering of the critical values is $S_1^5<T_1^5<S_2^5<S_3^5<T_2^5<T_3^5$.
}\label{n5}
\end{center}
\end{figure}

\begin{figure}[ht]
\begin{center}
\includegraphics[width=120mm,height=120mm]{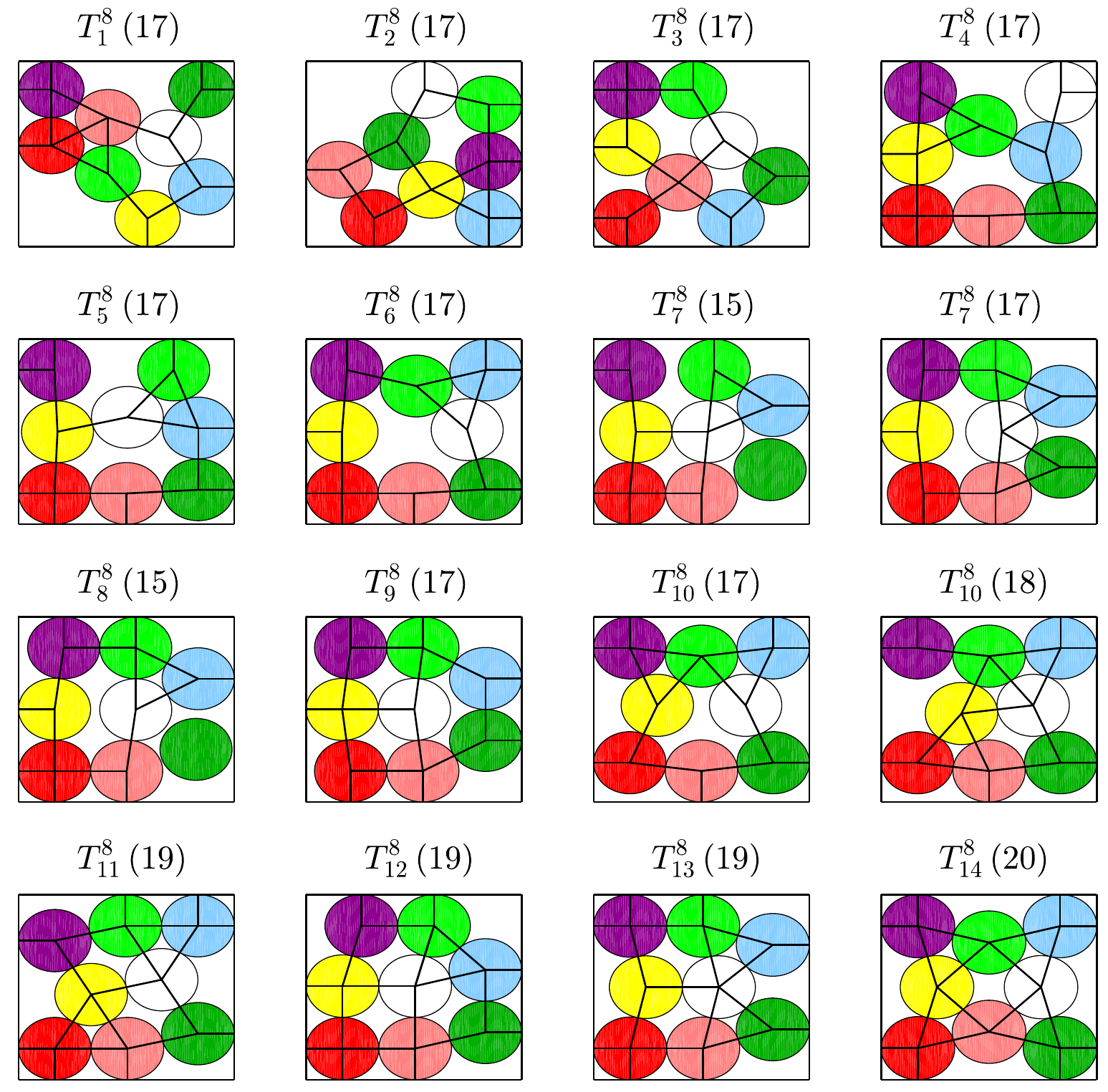}
\caption{Partially and fully jammed configurations of eight
disks with increasing radii from left to right and top to bottom. 
}\label{n8}
\end{center}
\end{figure}

The minimal example for multiple jammed configurations of different radii is
that of five disks in a square see Figure \ref{n5} (top row). In the first two
cases, the critical set consists of $4\cdot5!=480$ isolated points (one of the
four corners is distinguished). In the last case, the critical set consists of
$4!=24$ isolated points only, since this arrangement has a greater degree of
symmetry. We can  identify the cascade of radii at which the number of
connected components of $M_5^2(r)$ increases (Figure \ref{n5}, bottom row). We
know various  jammed configurations of six to eight disks (see Figures \ref{n6}, \ref{n7} and \ref{n8}).  We also find the minimal polynomials whose smallest positive roots are the maximal radii with the help of the \textsc{Mathematica} function \texttt{GroebnerBasis}, see Tables \ref{minpolys} and \ref{minpolyT58}.
 In cases where it is not obvious, a configuration can be tested for full jamming by Connelly's criterion \cite[Equation (1)]{Connelly2008}. Namely, if there is an infinitesimal motion $\x'=(\x_1',\,\dots\,,\x_n')$ of the configuration \mbox{$\x=(\x_1,\,\dots\,,\x_n)$} that satisfies
\begin{equation*}
(\x_i-\x_j)\cdot(\x_i'-\x_j')\ge0
\end{equation*}
for every pair $(i,j)$ of touching disks, then the motion should vanish, i.e.~$\x'=0$, for the configuration to be fully jammed. 

\begin{figure}[ht]
\begin{center}
\includegraphics[width=60mm,height=55mm]{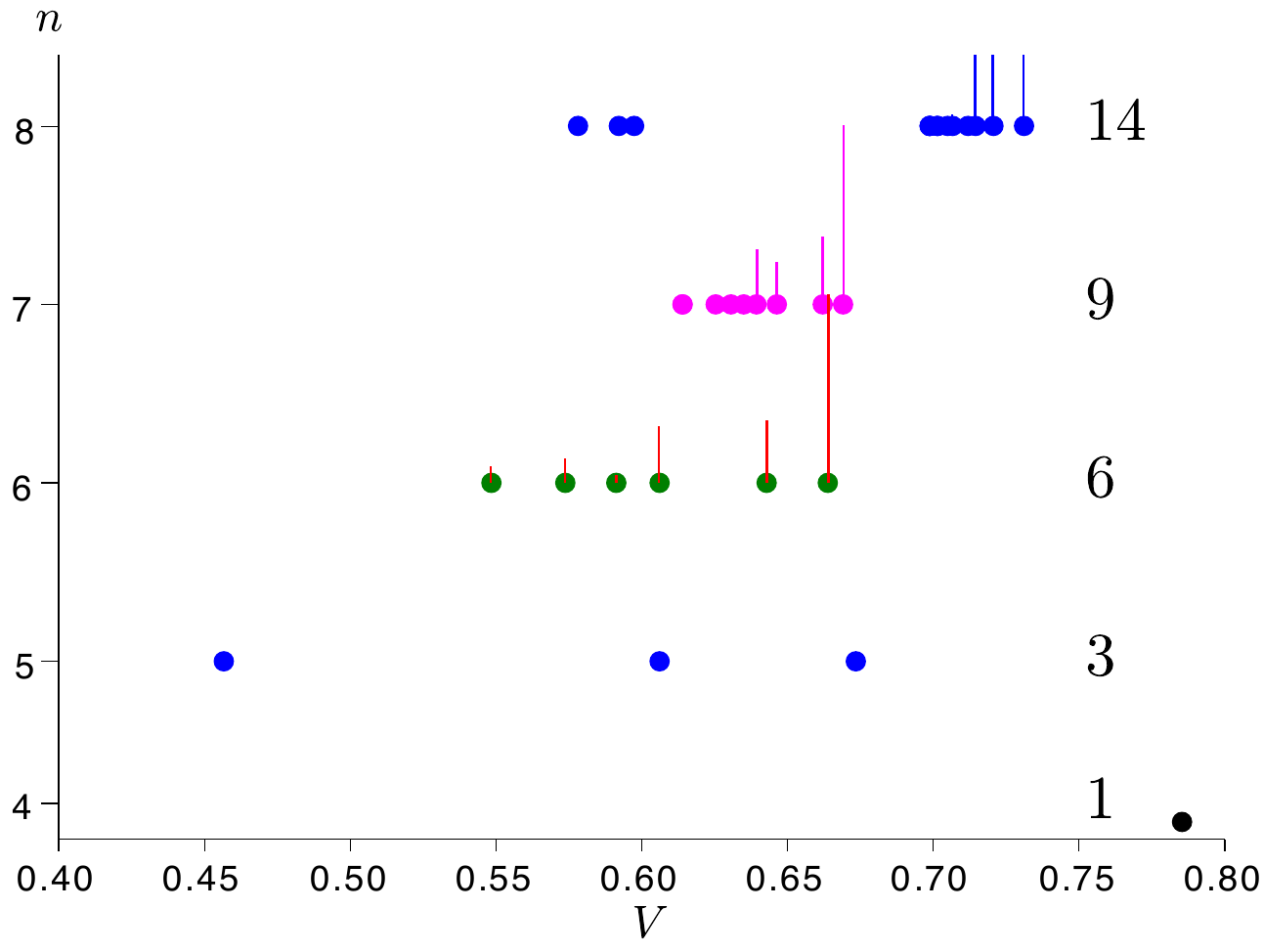}
\caption{The volume fractions realized by jammed packings of 4 to 8 disks. The histograms show the relative frequency of the terminal value of the maximization procedure after 5000 runs. }\label{trees_5_7}
\end{center}
\end{figure}

We repeat the maximization procedure described in Section \ref{Theorem} about $10^4$ times with initial choices $\x^0$ distributed uniformly in $[0,1]^{nd}$. Locally maximal configurations are compared to a list of known maxima (listed with all their images under symmetries of the square) by successive minimization of distances between the points of the new candidate and the known configuration.  Although the maximization procedure is not entirely deterministic due to the randomized symmetry breaking, this gives a rough measure of the relative sizes of the basins of attraction of the different maxima. It is a general observation that denser configurations are reached more frequently, see  the right panel of Figure \ref{trees_5_7}.

\begin{figure}[ht]
\begin{center}
\includegraphics[width=120mm]{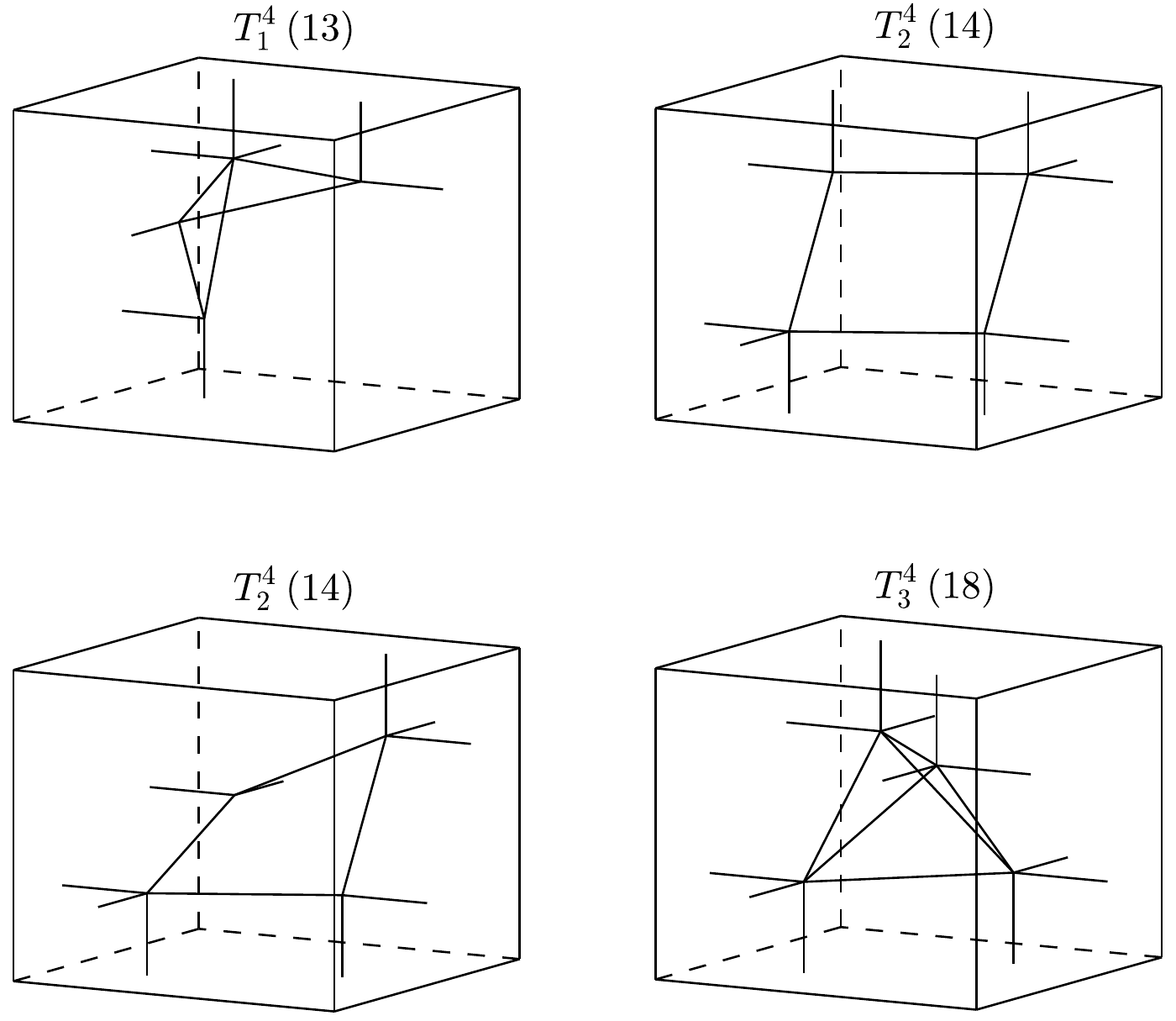}
\caption{Contact graphs of packings of four spheres in $[0,1]^3$ with increasing radii from left to right and top to bottom. Note that the two configurations labeled $T^4_2$ have identical contact graphs, but different arrangements of the spheres in the cube. }\label{n4d3}
\end{center}
\end{figure}

\begin{figure}[ht]
\begin{center}
\includegraphics[width=120mm]{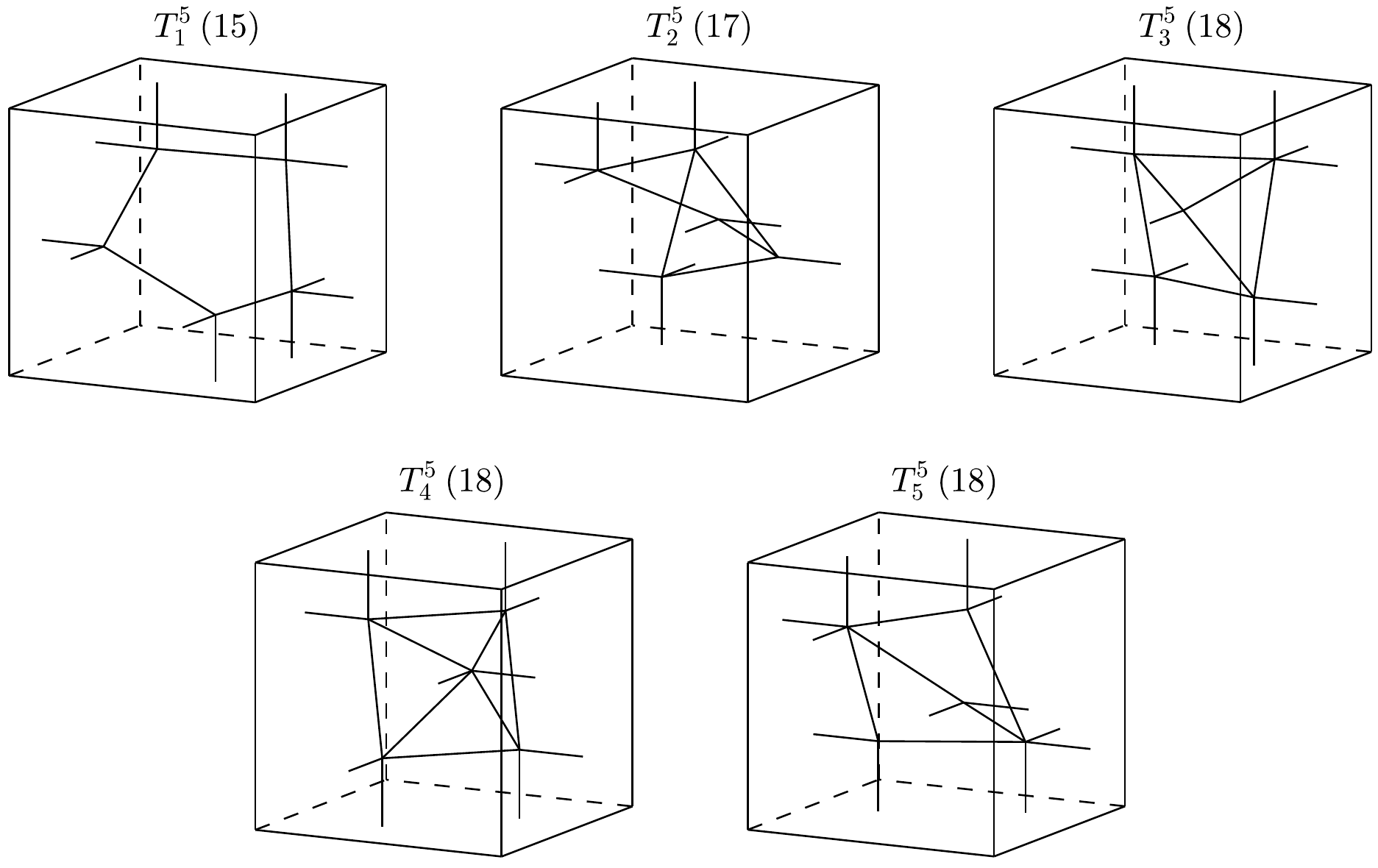}
\caption{Contact graphs of packings of five spheres in $[0,1]^3$ with increasing radii from left to right and top to bottom.  }\label{n5d3}
\end{center}
\end{figure}

\begin{figure}[ht]
\begin{center}
\includegraphics[width=120mm]{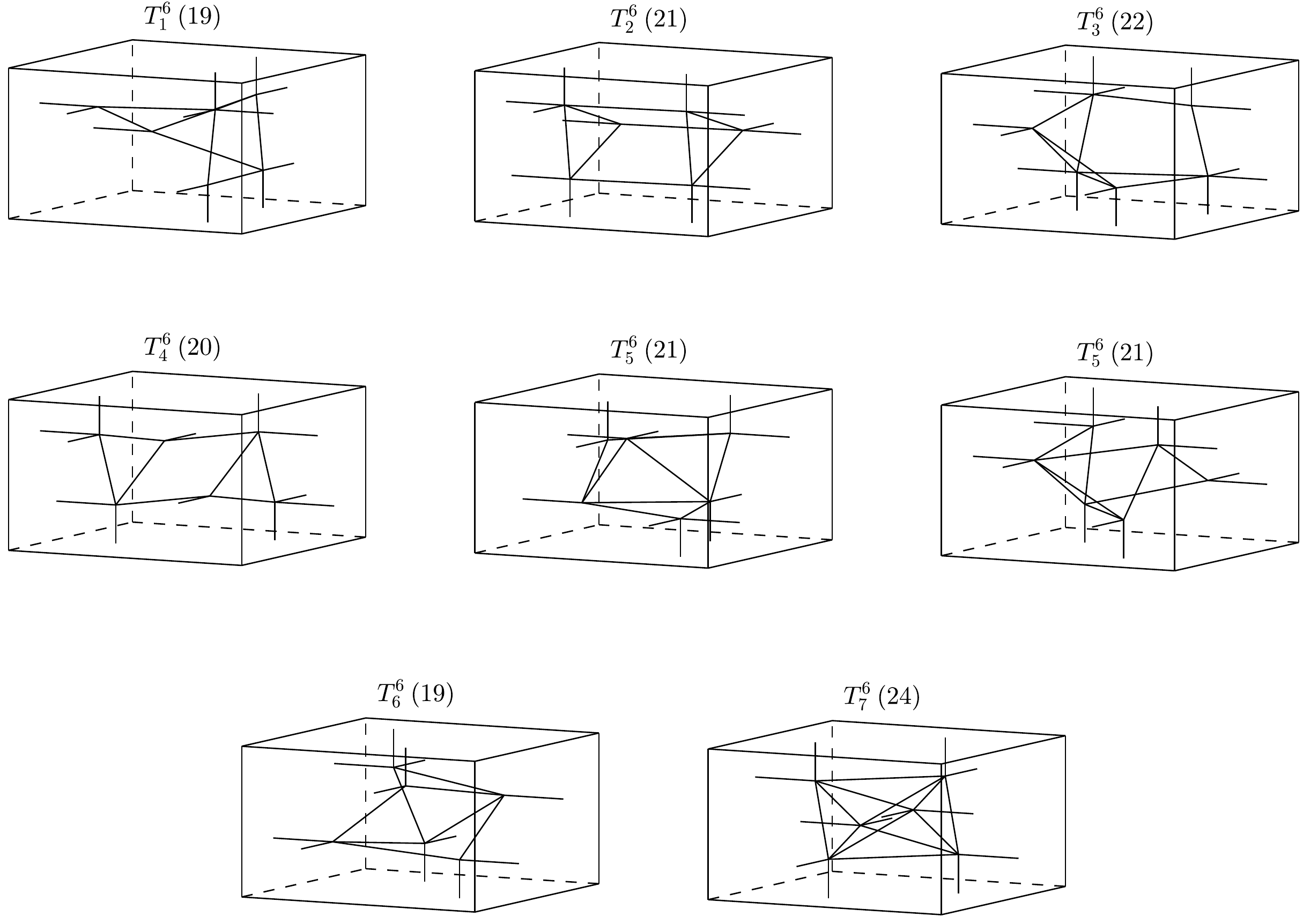}
\caption{Contact graphs of packings of six spheres in $[0,1]^3$ with increasing radii from left to right and top to bottom.  }\label{n6d3}
\end{center}
\end{figure}

Contact graphs of packings of four to six spheres in $[0,1]^3$ are given in Figures \ref{n4d3} - \ref{n6d3} and Tables  \ref{tabnd3} - \ref{tabnd3c}.

\begin{figure}[ht]
\begin{center}
\includegraphics[width=76mm,height = 36mm]{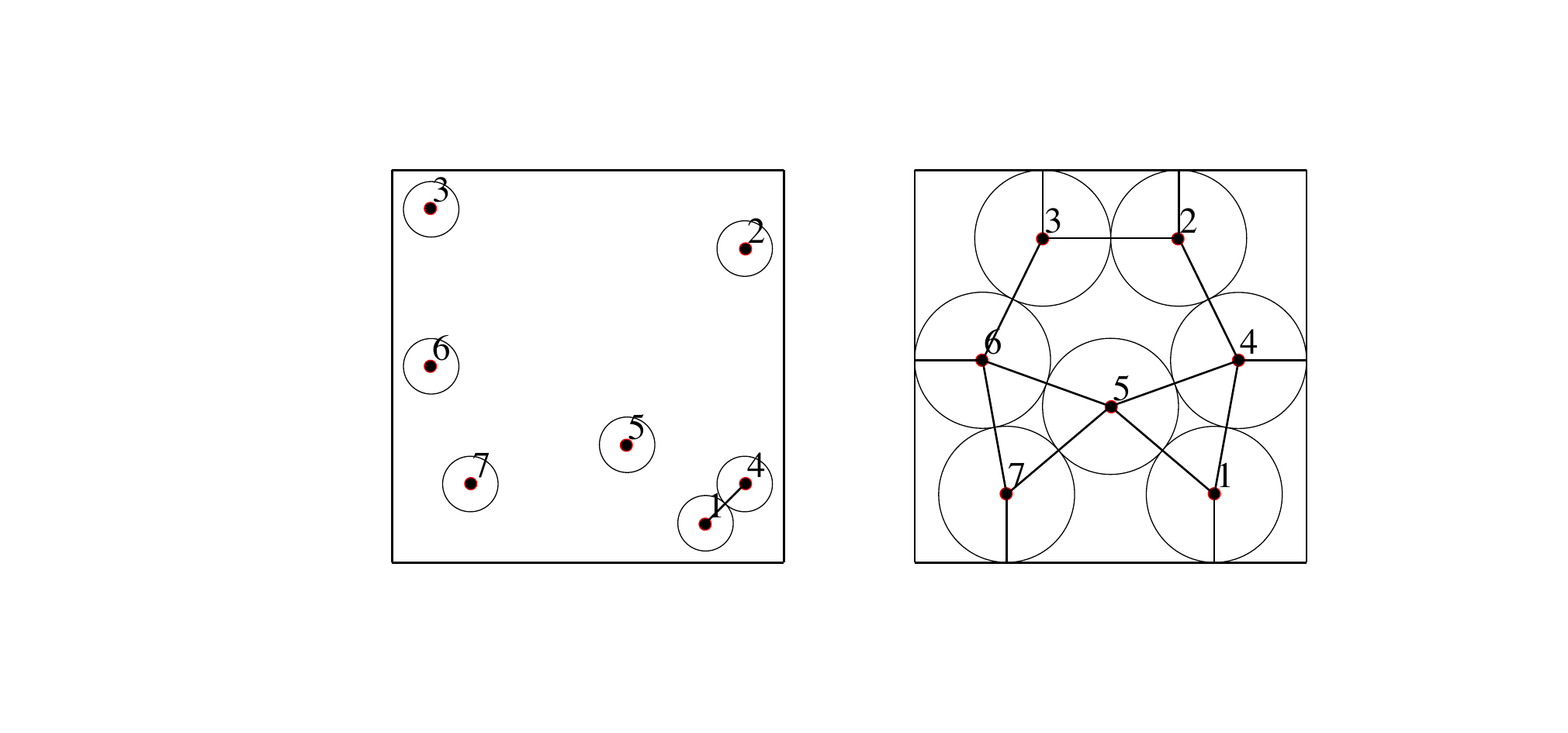}
\includegraphics[width=40mm,height = 36mm]{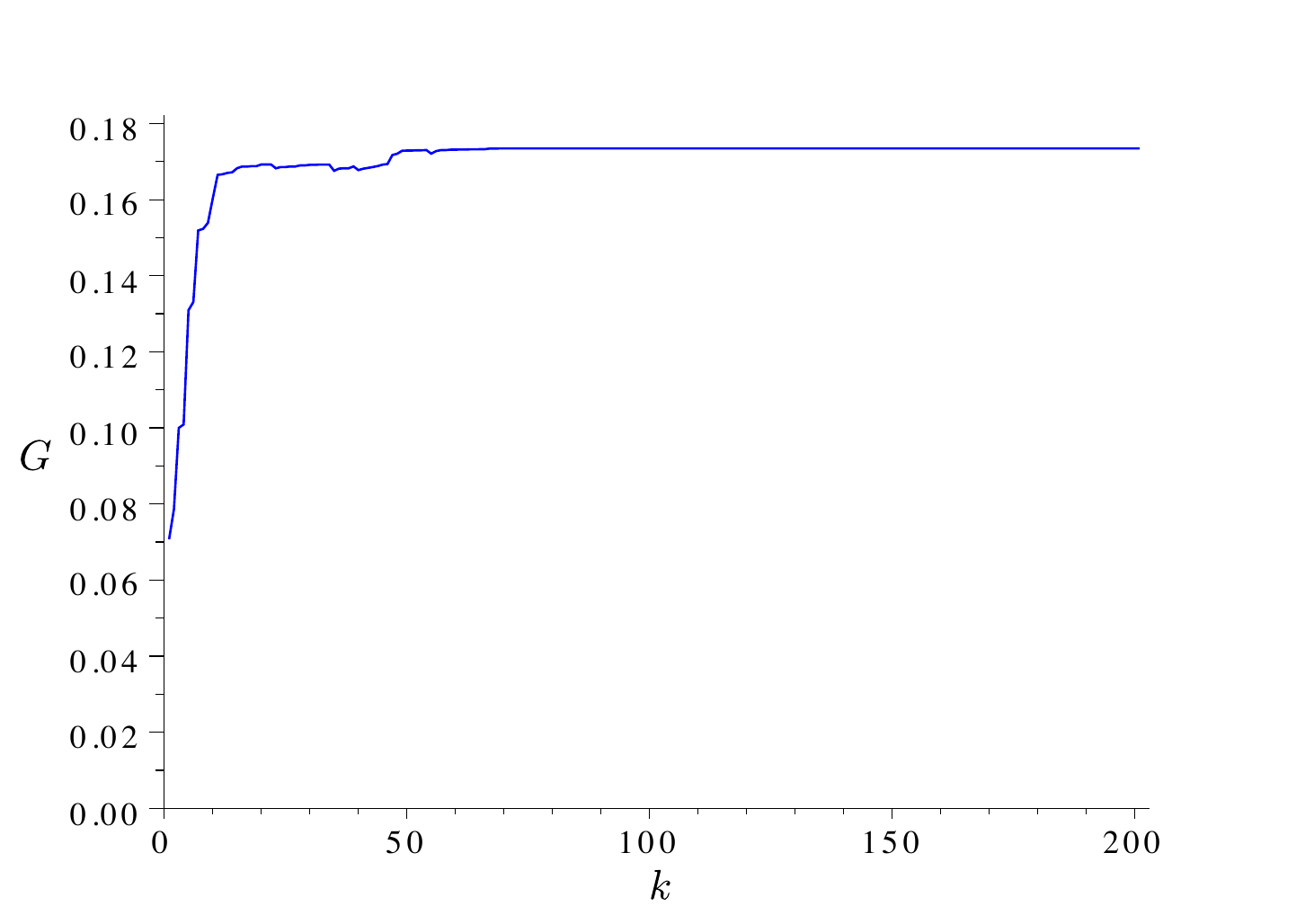}
\caption{{\it (Left)} An initial configuration of seven disks in the unit square and a local maximum provided by the search algorithm. {\it (Right)} The corresponding value of the function $G$ over 200 iterations of the line search procedure.  }\label{bench1}
\end{center}
\end{figure}

\begin{figure}[ht]
\begin{center}
\includegraphics[width=76mm,height = 36mm]{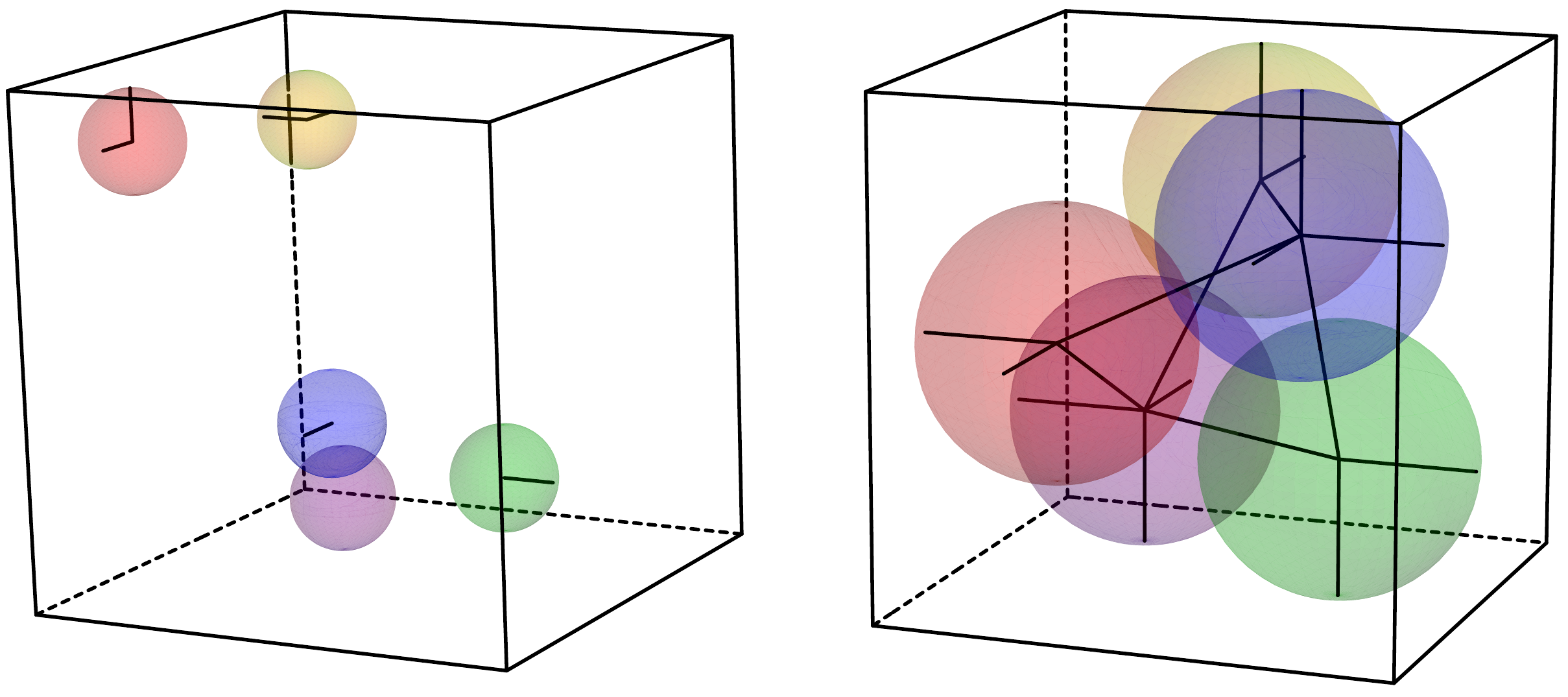}
\includegraphics[width=40mm,height = 36mm]{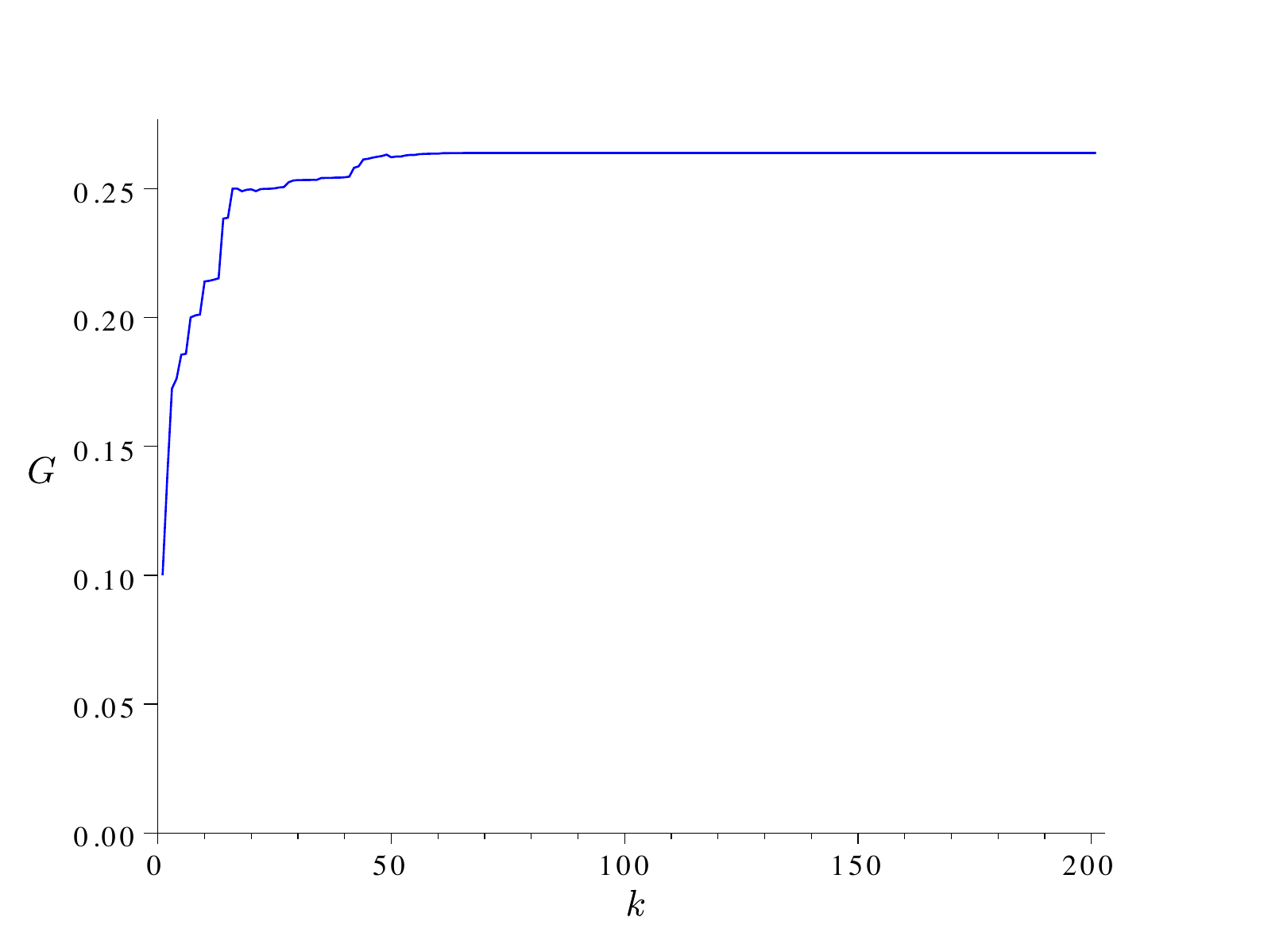}
\caption{{\it (Left)} An initial configuration of five spheres in the unit cube and a local maximum provided by the search algorithm. {\it (Right)} The corresponding value of the function $G$ over 200 iterations of the line search procedure.}\label{bench2}
\end{center}
\end{figure}

\begin{figure}[ht]
\begin{center}
\includegraphics[width=76mm,height = 36mm]{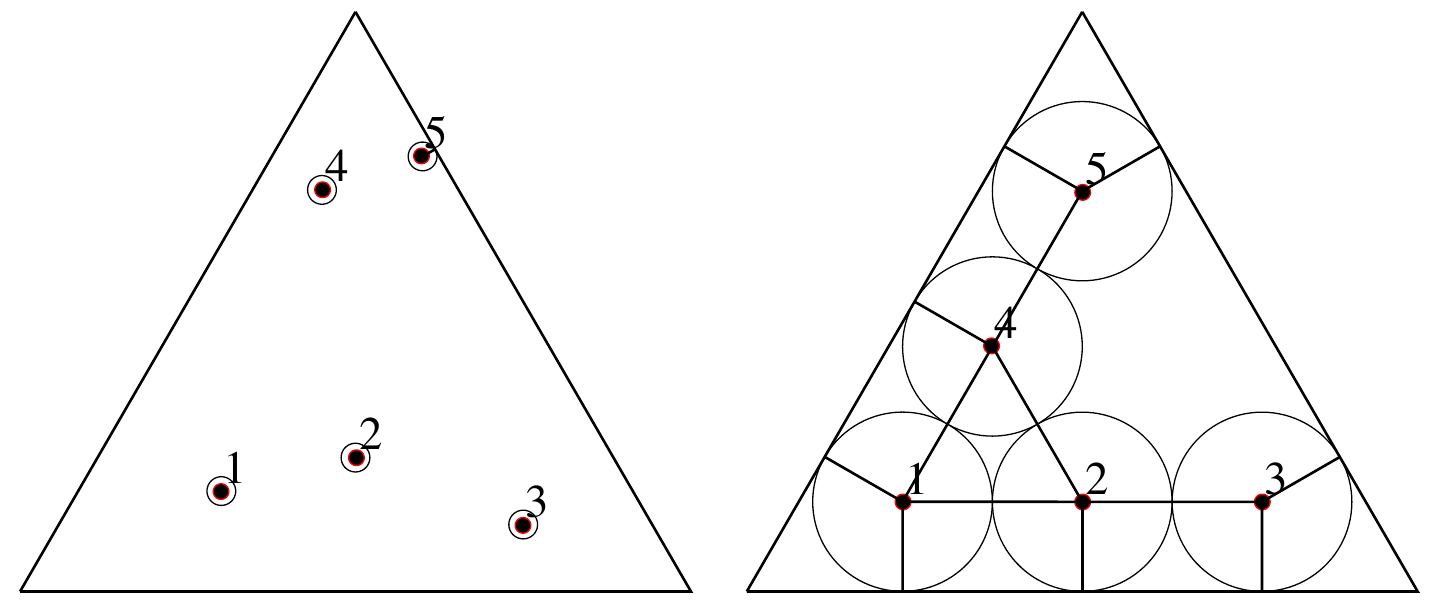}
\includegraphics[width=40mm,height = 36mm]{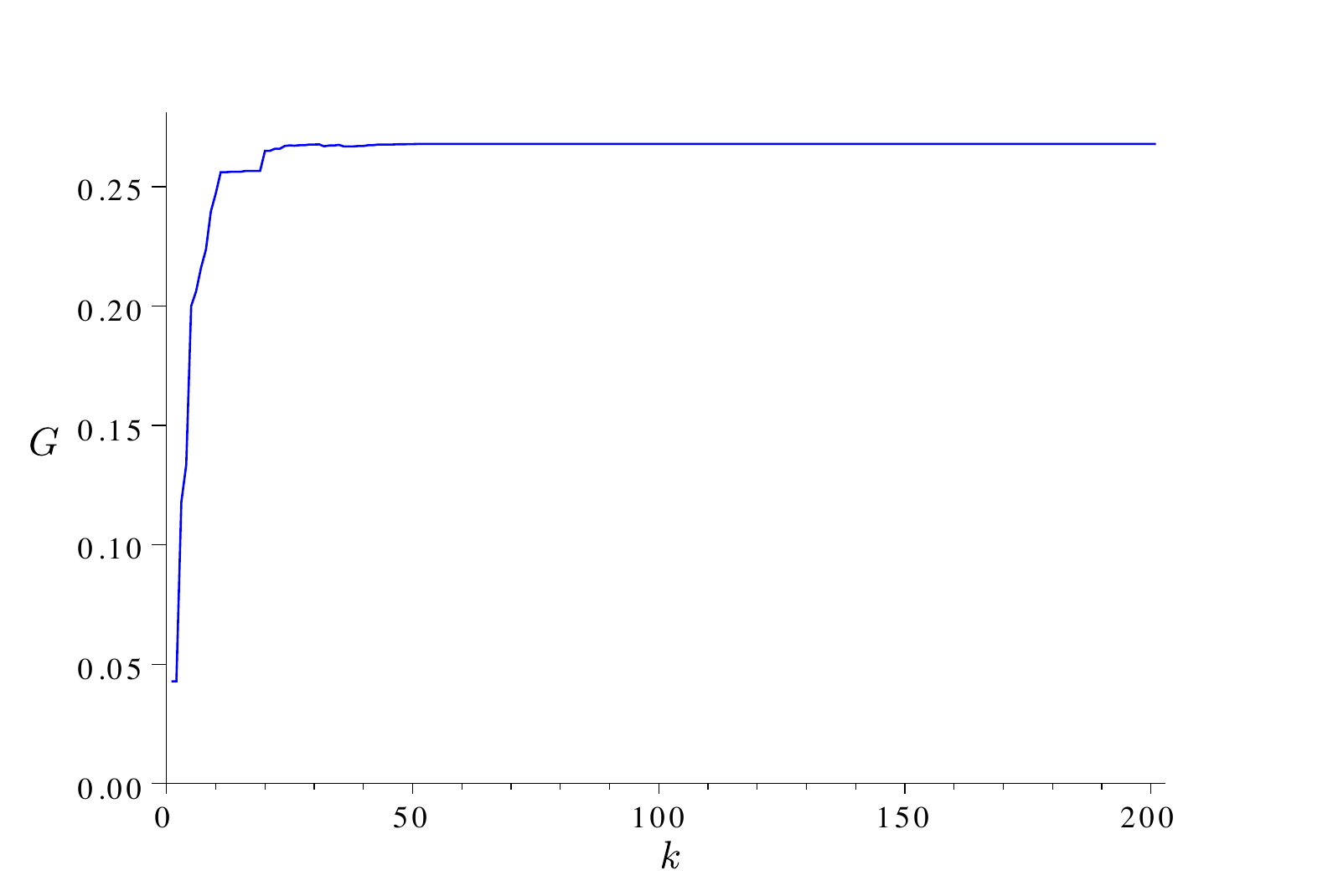}
\caption{{\it (Left)} An initial configuration of five disks in the an equilateral triangle of side lenght 2 and a local maximum provided by the search algorithm. {\it (Right)} The corresponding value of the function $G$ over 200 iterations of the line search procedure.}\label{bench3}
\end{center}
\end{figure}

We end this section with some test configurations that allow other researchers to benchmark their methods. It has to be stated again  that our algorithm has a random component and  the path of the algorithm to the  local maximum cannot be reproduced exactly. The first  is a configuration of seven disks in the square, the second is a configuration of five spheres in the cube, and the third is a configuration of five disks in an equilateral triangle of side length 2,
\begin{equation*}
\begin{aligned}
\x &= (0.8, 0.1; \,0.9, 0.8;\, 0.1, 0.9;\, 0.9, 0.2;\, 0.6, 0.3;\, 0.1, 0.5;\, 0.2, 0.2), \\
\x &= (0.1, 0.9, 0.8;\, 0.2, 0.1, 0.9;\, 0.9, 0.2, 0.3;\, 0.6, 0.1, 0.4;\, 0.5, 0.3, 0.2), \\
\x &= (0.6, 0.3;\, 1.0, 0.4;\, 1.5, 0.2;\, 0.9, 1.2;\, 1.2, 1.3).
\end{aligned}
\end{equation*}
Figures \ref{bench1} - \ref{bench3} show the results of runs with these initial configurations.  

\section{Discussion}

The above methodology may be applied to other geometrical setups, for example spheres in a regular simplex or spherical caps on the unit sphere $\mathbb{S}^d$ (also known as \textit{Tammes problem}). It is also possible to have multiple radius classes with fixed radius ratios. As stated in the introduction, the goal of the work is not to find the global maximum of the maximal radius function $G$, but rather as many local maxima as possible. Unfortunately, there is no way of knowing how many local maximal values there are, nor to know the ``basins of attraction'' of the local maxima for the given (or any other) search algorithm. There is a choice of alternative algorithms for a single  optimization of the function $G$, for example the gradient sampling algorithm \cite{Burke} or derivative-free methods \cite{Hare}. One also has to study alternative choices for the distribution of initial points (random or not), in order to find as many local maxima as possible. This touches upon the ``curse of dimensionality'' which is a  well-known problem in  numerical integration (for example over $[0,1]^d$), see \cite{Kuo} for a recent review.

It is possible to improve the performance of the search algorithm by exploiting the special structure of the objective function $G$ further. Among the $O(n^2)$ pairwise distances between sphere centers, many need not to be computed for finding the minimum. For example, one can subdivide the unit cube into $n$ cubes of equal side lengths and only compute those distances between sphere centers in adjacent subregions \cite{Lubachevsky1990}. This introduces some overhead, since membership lists of the subregions must be maintained. 

Although the data in Figures \ref{n5}, \ref{n6}, \ref{n7}, \ref{n8} and \ref{n5d3} are sparse, they allow to make some interesting observations. First we observe that no contact graph appears at more than one critical value. This leads to interesting combinatorial problems. Is it possible to give upper and lower bounds on the number of local maximal values of a minimum selection of convex functions as in Equation \eqref{lower_envelope}? A very coarse upper bound would be the number of subsets of the index set $\LL$ with $nd+1$ and more elements. Since the functions $\psi_{ik}$ and $\psi^{ik}$ are complementary to each other, at most one of these can be active and so a slightly reduced upper bound for minimally determined maxima would be
\begin{equation*}
\left(\begin{array}{c} n(n-1)/2+nd \\ nd+1 \end{array}\right).
\end{equation*}

A starting point of this work was the question which packing fractions are realizable by strictly jammed packings. For a jammed packing $\x^*$ of $n$ spheres in $[0,1]^d$, the packing fraction is 
\begin{equation*}
\phi_{n,d}=n\omega_d G(\x^*)^d, 
\end{equation*}
where $\omega_d$ is the volume of the $d$-dimensional unit ball.  Here only few general things can be said, apart from the numerical area spectra in the right panel of Figure \ref{trees_5_7}. In every dimension, there is at least one packing
fraction that is realized infinitely many times, namely $\frac{\omega_d}{2^d}$, which can be achieved by any number of $n^d$ spheres. 
Other repeated packing patterns are known for disks \cite{Lubachevsky1996,Nurmela97}. 

To end this paper with an open question, is there an open interval in $[0,1]$ of packing fractions that are realizable as limits of finite packings as $n\to\infty$?

\section*{Acknowledgments}
This work was partially supported by the NSF grant ``Collaborative Research: Predicting the Release Kinetics of Matrix Tablets'' (DMS 1016214). I thank Drs.~Boris Okun and Michael Hero (Department of Mathematical Sciences, University of Wisconsin - Milwaukee) and Michael Monagan  (Department of Mathematics, Simon Fraser University, Burnaby, British Columbia)  for helpful discussions.
I  thank the reviewer for the careful reading of the manuscript and the constructive remarks.

\bibliography{random_spheres}

\begin{thebibliography}{10}

\bibitem{DCDS_B}
B.~Baeumer, L.~Chatterjee, P.~Hinow, T.~Rades, A.~Radunskaya, and I.~Tucker.
\newblock Predicting the drug release kinetics of matrix tablets.
\newblock {\em Discr. Contin. Dyn. Sys. B}, {\bf 12}:261--277, 2009.
\newblock \href{http://arxiv.org/abs/0810.5323}{\texttt{arXiv:0810.5323}}.

\bibitem{Bernal1960a}
J.~D. Bernal.
\newblock Geometry of the structure of monatomic liquids.
\newblock {\em Nature}, {\bf 185}:68--70, 1960.

\bibitem{Bernal1960b}
J.~D. Bernal and J.~Mason.
\newblock Packing of spheres: co-ordination of randomly packed spheres.
\newblock {\em Nature}, {\bf 188}:910--911, 1960.

\bibitem{Boll}
D.~W. Boll, J.~Donovan, R.~L. Graham, and B.~D. Lubachevsky.
\newblock Improving dense packing of equal disk in a square.
\newblock {\em Electron. J. Combin.}, {\bf 7}:paper~\#46, 2000.

\bibitem{Burke}
J.~V. Burke, A.~S. Lewis, and M.~L. Overton.
\newblock A robust gradient sampling algorithm for nonsmooth, nonconvex
  optimization.
\newblock {\em SIAM J. Optim.}, {\bf 15}:751--779, 2005.

\bibitem{Carlsson}
G.~Carlsson, J.~Gorham, M.~Kahle, and J.~Mason.
\newblock Computational topology for configuration spaces of hard disks.
\newblock {\em Phys. Rev. E}, {\bf 85}:011303, 2012.
\newblock \href{http://arxiv.org/abs/1108.5719}{\texttt{arxiv:1108.5719}}.

\bibitem{Clarke1983}
F.~H. Clarke.
\newblock {\em Optimization and Nonsmooth Analysis}.
\newblock Wiley Interscience, New York, 1983.

\bibitem{Connelly2008}
R.~Connelly.
\newblock Rigidity of packings.
\newblock {\em European J. Combin.}, {\bf 29}:1862--1871, 2008.

\bibitem{Croft}
H.~T. Croft, K.~J. Falconer, and R.~K. Guy.
\newblock {\em Unsolved Problems in Geometry}.
\newblock Springer Verlag, New York, Berlin, Heidelberg, 1991.

\bibitem{Diaconis09}
P.~Diaconis.
\newblock The {Markov chain Monte Carlo} revolution.
\newblock {\em Bull. Amer. Math. Soc.}, {\bf 46}:179--205, 2009.

\bibitem{scilab}
{Digiteo Foundation, INRIA}.
\newblock \textsc{scilab}.
\newblock Available at \href{http://www.scilab.org}{\texttt{www.scilab.org}}.

\bibitem{Donev2004a}
A.~Donev, S.~Torquato, F.~H. Stillinger, and R.~Connelly.
\newblock Jamming in hard spheres and disk packings.
\newblock {\em J. Appl. Phys.}, {\bf 95}:989--999, 2004.

\bibitem{Gensane}
Th. Gensane.
\newblock Dense packings of equal spheres in a cube.
\newblock {\em Electron. J. Combin.}, {\bf 11}:paper~\#33, 2004.

\bibitem{Hardin}
D.~P. Hardin and E.~B. Saff.
\newblock Discretizing manifolds via minimum energy points.
\newblock {\em Notices Amer. Math. Soc.}, {\bf 51}:1186--1194, 2004.

\bibitem{Hare}
W.~L. Hare and M.~Macklem.
\newblock Derivative-free optimization methods for finite minimax problems.
\newblock {\em Optim. Methods Softw.}, {\it to appear}, 2011.

\bibitem{Joos}
A.~Jo\'os.
\newblock On the packing of fourteen congruent spheres in a cube.
\newblock {\em Geom. Dedicata}, {\bf 140}:49--80, 2009.

\bibitem{Kamien}
R.~D. Kamien and A.~J. Liu.
\newblock Why is random close packing reproducible?
\newblock {\em Phys. Rev. Lett.}, {\bf 99}:155501, 2007.
\newblock \href{http://arxiv.org/abs/cond-mat/0701343}{\texttt{
  arxiv:cond-mat/0701343}}.

\bibitem{Knott}
G.~M. Knott, T.~L. Jackson, and J.~Buckmaster.
\newblock Random packings of heterogeneous propellants.
\newblock {\em AIAA Journal}, {\bf 39}:678--686, 2001.

\bibitem{Kuo}
F.~Y. Kuo and I.~H. Sloan.
\newblock Lifting the curse of dimensionality.
\newblock {\em Notices Amer. Math. Soc.}, {\bf 52}:1320–--1329, 2005.

\bibitem{Loewen}
H.~L{\"o}wen.
\newblock Fun with hard spheres.
\newblock In {\em Statistical Physics and Spatial Statistics}, volume~{\bf 554}
  of {\em Lecture Notes in Physics}, pages 295--331, New York, Berlin,
  Heidelberg, 2000. Springer Verlag.

\bibitem{Lubachevsky1996}
B.~D. Lubachevsky and R.~L. Graham.
\newblock Repeated patterns of dense packings of equal disks in a square.
\newblock {\em Electron. J. Combin.}, {\bf 3}:paper~\#17, 1996.

\bibitem{Lubachevsky1990}
B.~D. Lubachevsky and F.~H. Stillinger.
\newblock Geometric properties of random disk packings.
\newblock {\em J. Stat. Phys.}, {\bf 60}:561--583, 1990.

\bibitem{Makela92}
M.~M. M{\"a}kela and P.~Neittaanm{\"a}ki.
\newblock {\em Nonsmooth Optimization}.
\newblock World Scientific, Singapore, 1992.

\bibitem{Maranas}
C.~D. Maranas, C.~Floudas, and P.~M. Pardalos.
\newblock New results in the packing of equal circles in a square.
\newblock {\em Discrete Math.}, {\bf 142}:287–--293, 1995.

\bibitem{Nurmela97}
K.~J. Nurmela and P.~R. J.~{\"O}sterg\aa rd.
\newblock Packing up to 50 equal circles in a square.
\newblock {\em Discrete Comput. Geom.}, {\bf 18}:111--120, 1997.

\bibitem{Nurmela99}
K.~J. Nurmela and P.~R. J.~{\"O}sterg\aa rd.
\newblock More optimal packing of equal circles in a square.
\newblock {\em Discrete Comput. Geom.}, {\bf 22}:439--457, 1999.

\bibitem{Peikert}
R.~Peikert, D.~W{\"u}rtz, M.~Monagan, and C.~de~Groot.
\newblock Packing circles in a square: a review and new results.
\newblock In {\em System Modelling and Optimization (Z{\"u}rich, 1991)},
  volume~{\bf 180} of {\em Lecture Notes in Control and Information Sciences},
  pages 45--54, Berlin, 1992. Springer Verlag.

\bibitem{Scott1960}
G.~D. Scott.
\newblock Packing of spheres: {P}acking of equal spheres.
\newblock {\em Nature}, {\bf 188}:908--909, 1960.

\bibitem{Torquato}
S.~Torquato, T.~M. Truskett, and P.~G. Debenedetti.
\newblock Is random close packing of spheres well defined?
\newblock {\em Phys. Rev. Lett.}, {\bf 84}:2064, 2000.
\newblock \href{http://arxiv.org/abs/cond-mat/0003416}{\texttt{
  arxiv:cond-mat/0003416}}.

\end{thebibliography}

\section*{Figures and Tables}


\begin{table}[ht]
\begin{center}
\begin{tabular}{|l|l|}
\hline
value & minimal polynomial\\
\hline
$T_1^5$ & $-1+8 t-8 t^2-32 t^3+32 t^4$ \\
$T_2^5$ & $1-12 t+52 t^2-100 t^3+73 t^4$ \\
$T_3^5$ & $-1-4 t+4 t^2$ \\
\hline
$T_1^6$ & $-1+8 t-8 t^2-32 t^3+32 t^4$ \\
$T_2^6$ & $1 - 8 t + 13 t^2$ \\
$T_3^6$ & $1-6 t+2 t^2$ \\
$T_4^6$ & $2-38 t+293 t^2-1180 t^3+2632 t^4-3104 t^5+1524 t^6$ \\
$T_5^6$ & $1-8 t+14 t^2$ \\
$T_6^6$ & $-13+52 t+92 t^2$ \\
\hline
$T_1^7$ & $37-472 t+1872 t^2-2368 t^3+832 t^4$ \\
$T_2^7$ & $4-64 t+368 t^2-896 t^3+769 t^4$ \\
$T_3^7$ & $1-16 t+92 t^2-224 t^3+193 t^4$ \\
$T_4^7$ & $13-112 t+208 t^2$ \\
$T_5^7$ & $25-440 t+4656 t^2-34176 t^3+144896 t^4-311296 t^5+262144 t^6$ \\
$T_6^7$ & $-1+8 t-8 t^2-32 t^3+32 t^4$ \\
$T_7^7$ & $5-138 t+1650 t^2-11280 t^3+49316 t^4-147424 t^5+316132 t^6$ \\
               & $-495168 t^7+547136 t^8-400512 t^9+196096 t^{10}-58368 t^{11}+9216 t^{12}$ \\
$T_8^7$ & $481-5680 t+22856 t^2-40416 t^3+32016 t^4-9344 t^5+256 t^6$ \\
$T_9^7$ & $1-8 t+13 t^2$ \\
\hline
$T_1^8$ & see Table \ref{minpolyT58} \\
$T_2^8$ & $9-342 t+5850 t^2-59538 t^3+400966 t^4-1876622 t^5+6223518 t^6$ \\
               & $-14613362 t^7+23977405 t^8-27600528 t^9+25173396 t^{10}$ \\
               & $-22602672 t^{11}+14091408 t^{12}$ \\
$T_3^8$ & $20-300 t+1621 t^2-4152 t^3+5832 t^4-4320 t^5+1296 t^6$ \\
$T_4^8$ & $1-50 t+1134 t^2-15474 t^3+142078 t^4-929874 t^5+4481482 t^6$ \\
               & $-16190694 t^7+44168181 t^8-90780660 t^9+138781596 t^{10}$ \\
               & $-153497992 t^{11}+116495620 t^{12}-54495056 t^{13}+11899792 t^{14}$ \\
$T_5^8$ & see Table \ref{minpolyT58} \\
$T_6^8$ & $1-24 t+248 t^2-1472 t^3+5640 t^4-14624 t^5+25248 t^6$ \\ &$-26240 t^7+12304 t^8$ \\
$T_7^8$ & $37-472 t+1872 t^2-2368 t^3+832 t^4$ \\
$T_8^8$ & $4-64 t+368 t^2-896 t^3+769 t^4$ \\
$T_9^8$ & $2-30 t+155 t^2-310 t^3+169 t^4$ \\
$T_{10}^8$ & $-73+584 t-344 t^2-3296 t^3+368 t^4$ \\
$T_{11}^8$ & $1-12 t+52 t^2-120 t^3+148 t^4$ \\
$T_{12}^8$ & $1-16 t+92 t^2-224 t^3+193 t^4$ \\
$T_{13}^8$ & $13-112 t+208 t^2$ \\
$T_{14}^8$ & $-1+8 t-8 t^2-32 t^3+32 t^4$ \\
\hline
\end{tabular}
\caption{Minimal polynomials of the locally maximal radii of disks in the unit square, which is always the first positive root.}\label{minpolys}
\end{center}
\end{table}

\begin{landscape}
\begin{table}[ht]
\begin{center}
\begin{tabular}{|l|l|}
\hline
value & minimal polynomial\\
\hline
$T_1^8$ & $4 - 232 t + 6400 t^2 - 111960 t^3 + 1397292 t^4 - 13270344 t^5 +99892240 t^6 - 612729996 t^7 + 3122576292 t^-13397230544 t^9$ \\
& $ + 48791818124 t^{10} - 151444582764 t^{11} + 400899869393 t^{12}-903644665080 t^{13}+1729356481328 t^{14}-2800032538112 t^{15} $ \\
& $+3819006125024 t^{16} - 4359110782976 t^{17} + 4115559413504 t^{18} -3146327596032 t^{19}+ 1874189155584 t^{20}$ \\
& $- 805545011200 t^{21} +203203010560 t^{22} - 5307301888 t^{23} - 4813930496 t^{24} - 
 7357071360 t^{25} + 3466854400 t^{26}$ \\ 
\hline
$T_5^8$ & $547600-54665280 t+2605016436 t^2-78880971596 t^3+1704389979277 t^4-27976098788988 t^5+362612144032096 t^6$ \\ 
& $-3808672735525104 t^7+33009962234551198 t^8-239162479838195956 t^9+1462102546014264124 t^{10}$ \\
& $-7592277950931300452 t^{11}+33635593700920570705 t^{12}-127460395807424738280 t^{13}+413510439834837377472 t^{14}$ \\ 
& $-1147656676176223600832t^{15}+2718317154671720808288 t^{16}-5471314759394629231104 t^{17}$ \\ 
& $+9297731342945147042048 t^{18}-13217748559820732803072 t^{19}+15518739633438654030080 t^{20}$ \\ 
& $-14779772068382475778048 t^{21}+11127757838288057618432 t^{22}-6372552344848630513664 t^{23}$ \\
& $+2607643199804769697792 t^{24}-679025777716818345984 t^{25}+84539217262821769216 t^{26}$\\
\hline
\end{tabular}
\caption{Minimal polynomials of two configurations of 8 disks in the unit square.}\label{minpolyT58}
\end{center}
\end{table}
\end{landscape}
\begin{table}[ht]
\begin{center}
\begin{tabular}{|l|l|}
\hline
value & minimal polynomial\\
\hline
$T_1^4$ & $1-16 t+104 t^2-352 t^3+704 t^4-1024 t^5+1216 t^6-768 t^7+64 t^8$ \\
$T_2^4$ & $5-20 t+4 t^2$ \\
$T_3^4$ & $5-4 t+2 t^2$ \\
\hline
$T_1^5$ & $13-104 t+240 t^2-128 t^3+t^4$ \\
$T_2^5$ & $1-8 t+20 t^2-16 t^3+16 t^6$ \\
$T_3^5$ & $1-16 t+104 t^2-352 t^3+704 t^4-1024 t^5+1216 t^6-768 t^7+64 t^8$ \\
$T_4^5$ & $9-36 t+4 t^2$ \\
$T_5^5$ & $5-20 t+4 t^2$ \\
\hline
$T_1^6$ & $628849-35215544 t+920916520 t^2-14925909856 t^3$ \\
& $+167788673872 t^4-1387008330496 t^5+8720616629536 t^6$ \\
& $-42537931528576 t^7+162698001135232 t^8-489697674935296 t^9$ \\
& $+1156878795748352 t^{10}-2126444318666752 t^{11}$ \\
& $+2994825062826240 t^{12}-3163357034848256 t^{13}$ \\
& $+2454984222679040 t^{14}-1420158379393024 t^{15}$ \\
& $+698664367800320 t^{16}-324380948299776 t^{17}$ \\
& $+55492564615168 t^{18}+84665382207488 t^{19}$ \\
& $-46125516718080 t^{20}-14650005520384 t^{21}$ \\
& $+9630847598592 t^{22}+3387487682560 t^{23}-753326358528 t^{24}$ \\
& $-527777660928 t^{25}-87879057408 t^{26}-4697620480 t^{27}$ \\
& $+16777216 t^{28}$\\
$T_2^6$ & $13-104 t+240 t^2-128 t^3+t^4$ \\
$T_3^6$ & $1-8 t+20 t^2-16 t^3+t^4$ \\
$T_4^6$ & $1-16 t+96 t^2-256 t^3+544 t^4-2304 t^5+5632 t^6-4096 t^7+256 t^8$ \\
$T_5^6$ & see Table \ref{tabnd3b} \\
$T_6^6$ & see Table \ref{tabnd3c}  \\
$T_7^6$ & $9-36 t+4 t^2$ \\
\hline
\end{tabular}
\caption{Minimal polynomials of the locally maximal radii of four to six spheres in three dimensions.}\label{tabnd3}
\end{center}
\end{table}
\begin{landscape}
\begin{table}[ht]
\begin{center}
\begin{tabular}{|l|l|}
\hline
value & minimal polynomial\\
\hline
$T_5^6$ & $43046721-6198727824 t+630697172688 t^2-47220701129280 t^3+2869181328420576 t^4-146191375492879104 t^5$ \\
               & $+6372404443192839552 t^6-240984441891113614848 t^7+7961701503621774723840 t^8$ \\
              & $-230636044887016540962816 t^9+5866602757184413918016512 t^{10}-130932539887133984695988224 t^{11}$ \\
              & $+2559019173189523030256331776 t^{12}-43734819932165772392688885760 t^{13}+653541214060594313354395664384 t^{14}$ \\
              & $-8550436965155206662510296694784 t^{15}+98165697119101462453344534249472 t^{16}$ \\
              & $-991571274153891901202702482538496 t^{17}+8833485388097631073100752556359680 t^{18}$ \\
              & $-69525727421221420895499959892836352 t^{19}+483764673966036336113525690731593728 t^{20}$ \\
             & $-2972516284985985180613114101432844288 t^{21}+16069372775596251285701792968987639808 t^{22}$ \\
             & $-75812442485694868130653849894440140800 t^{23}+307081175567122027628500848071811661824 t^{24}$ \\
             & $-1031282570129898206591875362924588433408 t^{25}+2622189007453528772245078823295934005248 t^{26}$ \\
             & $-3338923015726622930149069530627280207872 t^{27}-11012791793216708754824602984394761699328 t^{28}$ \\
             & $+97818451430492140307456629592329077915648 t^{29}-402696023773161755120692097679483404812288 t^{30}$ \\
             & $+1046044618270399290113987926669989645385728 t^{31}-980215765887707705641630145603600574840832 t^{32}$ \\
             & $-7119128106689314956495645921038577656922112 t^{33}+53441256937086617576117288148322014126407680 t^{34}$ \\
             & $-233824996892750458373335403636306331966636032 t^{35}+799085060891855426598810257272049763126083584 t^{36}$ \\
            & $-2303121259038575398192072066905831569715363840 t^{37}+5782966042491046620994328152688518944929087488 t^{38}$ \\
            & $-12871172818830139108473860525875156291896213504 t^{39}+25681740870347460045615987769224908820935868416 t^{40}$ \\
            & $-46377196420229066828149274524764670932992131072 t^{41}+76619508043582810233618197399453313671742095360 t^{42}$ \\
            & $-117537858136343126076716427902931103070472372224 t^{43}+170967916039096605933568209795505191939408396288 t^{44}$ \\
            & $-241931155777388500939124812954240097081893060608 t^{45}+340579887823077792142831763033997285799424425984 t^{46}$ \\
            & $-480355423296991574111819017307414871602785943552 t^{47}+670021674127675127013520018450556330414717272064 t^{48}$ \\
            & $-901581261456833478662343962164560901262740029440 t^{49}+1141580027710048532150300585124636125113316540416 t^{50}$ \\
            & $-1334829908839129036976169221379519688531163217920 t^{51}+1423481041866669106473377321025283045342375313408 t^{52}$ \\
           & $-1373330838629924939210784174249278223105556217856 t^{53}+1191863863383147017298364999858742642271106105344 t^{54}$ \\
           & $-926092087997101376832544535909299980694808690688 t^{55}+641234349803120868990847434687858297522033262592 t^{56}$ \\
           & $-393556579169074746491083478956830537233052205056 t^{57}+212731128094261011021317276200910436236989038592 t^{58}$ \\
\end{tabular}
\end{center}
\end{table}

\begin{table}[ht]
\begin{center}
\begin{tabular}{|l|l|}
\phantom{$T_5^6$}  
& $-100458492661460075376734095675283648912658792448 t^{59}+41019777270104513941100024681580239760547381248 t^{60}$ \\
& $-14289060656450598313079093887849494861150420992 t^{61}+4170552271286383771010770917305457563117027328 t^{62}$ \\
& $-994760408802755649711885615291786528284475392 t^{63}+186970979583422823433127549415424878807351296 t^{64}$ \\
& $-26161124413996589591885732626730180319117312 t^{65}+2474139291394158214389081922316174140375040 t^{66}$ \\
& $-136063865643192056740114801369454796079104 t^{67}+7040053447185167255747299754968998739968 t^{68}$ \\
& $-1668157208055489578064673415280168271872 t^{69}+271460796970985814401428444558756151296 t^{70}$ \\
& $-19893277926682173557147898049229488128 t^{71}+541847633547106597831607476587331584 t^{72}$ \\
\hline
\end{tabular}
\caption{Minimal polynomial of configuration $T_5^6$.}\label{tabnd3b}
\end{center}
\end{table}
\end{landscape}

\begin{table}[ht]
\begin{center}
\begin{tabular}{|l|l|}
\hline
value & minimal polynomial\\
\hline
$T_6^6$ 
& $17850625-3427320000 t+327817998528 t^2-20821826923264 t^3$ \\
& $+987726417983232 t^4-37313693603198976 t^5$ \\
& $+1168907639294074880 t^6-31219963077836095488 t^7$ \\
& $+725429602195873475072 t^8-14890754410484706631680 t^9$ \\
& $+273272074409891243776000 t^{10}-4526854896736132756807680 t^{11}$ \\
& $+68222063007549667368159232 t^{12}$ \\
& $-941472366257163248946987008 t^{13}$ \\
& $+11962606417179071943721172992 t^{14}$ \\
& $-140607566392008233399332569088 t^{15}$ \\
& $+1534983670264170814477055475712 t^{16}$ \\
& $-15618180998282167274579315654656 t^{17}$ \\
& $+148566699694213641005963068178432 t^{18}$ \\
& $-1324815392234110934918500019863552 t^{19}$ \\
& $+11101518878460160647263724506316800 t^{20}$ \\
& $-87608344468411144802040068910874624 t^{21}$ \\
& $+652374475492830725841694709689876480 t^{22}$ \\
& $-4592111614294706899389248431215083520 t^{23}$ \\
& $+30605659141869032519917088312592695296 t^{24}$ \\
& $-193427325865007346716195666391281434624 t^{25}$ \\
& $+1160816883779208050449949294650825113600 t^{26}$ \\
& $-6623654558628870509077400504674421833728 t^{27}$ \\
& $+35977991812578058500230628908181005271040 t^{28}$ \\
& $-186234401198438304242849849705495082827776 t^{29}$ \\
& $+919617920841624229559217715864229182963712 t^{30}$ \\
& $-4335933763403022506485823748571063052664832 t^{31}$ \\
& $+19536660186433782946138183434741460248821760 t^{32}$ \\ 
& $-84185159953565880143925521370127290339426304 t^{33}$ \\ 
& $+347153909291483094403651698938794743133372416 t^{34}$ \\
& $-1370726446384907196824324909361586404962861056 t^{35}$ \\
& $+5184639756910668175364231224718241744740679680 t^{36}$ \\
& $-18792199494194327747086885128090261375386910720 t^{37}$ \\
& $+65288329365056041796420742462877146325874376704 t^{38}$ \\
& $-217448630736569356568147184944646006522920304640 t^{39}$ \\
& $+694326625113846229591979032844468100187769798656 t^{40}$ \\
& $-2125402526658739058303714747630549925702678872064 t^{41}$ \\
& $+6236455238395122619884513811694349682546780930048 t^{42}$ \\
& $-17537647049426855545851198344255050069316202397696 t^{43}$ \\
& $+47253288368880129538720061055878544735369533849600 t^{44}$ \\
& $-121951401628459296647005792037531770690788067901440 t^{45}$ \\
& $+301357940570851601462385042150371731014485373091840 t^{46}$ \\
& $-712770888641580130614333642358451074651177636331520 t^{47}$ \\
& $+1612891265948281610269936395233794346617846429122560 t^{48}$ \\
& $-3490216403945726551834895659369234775060744696233984 t^{49}$ \\
& $+7219138516017948554073970044709643998449131791056896 t^{50}$ \\
& $-14265527844268055655822641834843546109671242475765760 t^{51}$ \\
& $+26917438040439551303328998818949515655121390821965824 t^{52}$ \\
\end{tabular}
\end{center}
\end{table}

\begin{table}[ht]
\begin{center}
\begin{tabular}{|l|l|}
\phantom{$T_6^6$}  
& $-48471547386879781831770952534709393150828442619478016 t^{53}$ \\
& $+83253087740653914752431225193771295139744614707101696 t^{54}$ \\
& $-136306032151071672536241697958629229037463749324177408 t^{55}$ \\
& $+212597617702783995411019367614122565664623794996641792 t^{56}$ \\
& $-315676397473038371399275005328607494391181053884628992 t^{57}$ \\
& $+445923436728234832244994572550328731272800693297086464 t^{58}$ \\
& $-598806842242809930669704195963713773883830547081330688 t^{59}$ \\
& $+763786818077750680401448080069900372795670519535894528 t^{60}$ \\
& $-924572288420772748644209019138331641051138898977095680 t^{61}$ \\
& $+1061175825386612849942695930583393999857612398744567808 t^{62}$ \\
& $-1153644758729836865896398564033754917894620519247904768 t^{63}$ \\
& $+1186634969364566686182010652820942165674057471044354048 t^{64}$ \\
& $-1153455507584389481769581120122977784941191436716474368 t^{65}$ \\
& $+1058161638153990547563625448826876248032168945639227392 t^{66}$ \\
& $-914835663828729166165471137790642289331085778821840896 t^{67}$ \\
& $+744187348818469308993553344010646859561311884061704192 t^{68}$ \\
& $-568596753514594051724592403539423821769311269813223424 t^{69}$ \\
& $+407248835742964482814890145907615762308072363154997248 t^{70}$ \\
& $-272836210124234733862768949209622314217163293068361728 t^{71}$ \\
& $+170557758832305386465277389754206931191924956380266496 t^{72}$ \\
& $-99215033862300201842972147687337951993209296674357248 t^{73}$ \\
& $+53539566761294723788199878504166278615696624851091456 t^{74}$ \\
& $-26707644285492858965320763857969512827073056556974080 t^{75}$ \\
& $+12266288629855639900270649057381011218333129081094144 t^{76}$ \\
& $-5162943193820983687414910787305081452370186372907008 t^{77}$ \\
& $+1980868601852602989513421247035536681150408206319616 t^{78}$ \\
& $-688426326927744918898952117680288113777061687984128 t^{79}$ \\
& $+215116218168063036670267319908135566372063166857216 t^{80}$ \\
& $-59901796639027703320498208900405319734826239000576 t^{81}$ \\
& $+14705581801215668409771970575321818506353256169472 t^{82}$ \\
& $-3141086480295116234155128029356961578965813690368 t^{83}$ \\
& $+574378710574108291290991339995528442291169525760 t^{84}$ \\
& $-88169677048730063219192101145477751204231512064 t^{85}$ \\
& $+11117813722503308733184750051934924875071750144 t^{86}$ \\
& $-1135155931044179239223315832239384127376195584 t^{87}$ \\
& $+97001378040253342688399202622900447232393216 t^{88}$ \\
& $-8255498437868897949003417002125119706890240 t^{89}$ \\
& $+870913736231297124613164904984002472968192 t^{90}$ \\
& $-98453434885687096630438107630143209996288 t^{91}$ \\
& $+9179886398371767408215264742671120334848 t^{92}$ \\
& $-651588097346297563395669662197706915840 t^{93}$ \\
& $+35581606297548954230865050104772952064 t^{94}$ \\
& $-1423200492877425383193585529990938624 t^{95}$ \\
& $+29623222488105672354653249885700096 t^{96}$ \\
\hline
\end{tabular}
\caption{Minimal polynomial of configuration $T_6^6$.}\label{tabnd3c}
\end{center}
\end{table}

\end{document}